\documentclass[10pt]{article}

\input diagram.tex
\usepackage{amsmath, amscd} 
\usepackage[mathscr]{eucal}
\usepackage{amssymb} 
\usepackage{theorem}
\usepackage{xspace}

\theoremstyle{change}  
\newtheorem{theorem}{Theorem}[section] 
\newtheorem{lemma}[theorem]{Lemma}  
\newtheorem{proposition}[theorem]{Proposition}
\newtheorem{corollary}[theorem]{Corollary}
\theorembodyfont{\rmfamily}  
\newtheorem{remark}[theorem]{Remark}
\newtheorem{example}[theorem]{Example}

\newtheorem{definition}[theorem]{Definition}
\newtheorem{notation}[theorem]{Notation}
\newtheorem{nothing}[theorem]{} 

\newenvironment{proof}{\noindent{\bf Proof}\ }{\qed\bigskip}

\renewcommand{\le}{\leqslant} 

\newcommand{\calC}{\mathcal{C}}
 
\newcommand{\calD}{\mathcal{D}}                      
\newcommand{\calE}{\mathcal{E}}
\newcommand{\calF}{\mathcal{F}}

\newcommand{\calG}{\mathcal{G}}

\newcommand{\calS}{\mathcal{S}}

\newcommand{\catfont}{\mathsf}

\newcommand{\CC}{\mathbb{C}}
\newcommand{\CCF}{\lexp{\CC}{F}}

\newcommand{\cdotH}{\mathop{\cdot}\limits_{H}}

\newcommand{\ConfR}{\mathrm{Con}_R^f}
\newcommand{\ConfmuR}{\mathrm{Con}_R^{f,\mu}}

\newcommand{\defl}{\mathrm{def}}

\newcommand{\Dphi}{D^{\varphi}}

\newcommand{\Hom}{\mathrm{Hom}}

\newcommand{\id}{\mathrm{id}}

\newcommand{\ind}{\mathrm{ind}}

\newcommand{\infl}{\mathrm{inf}}

\newcommand{\isog}{\mathrm{iso}}

\newcommand{\lexp}[2]{\setbox0=\hbox{$#2$} \setbox1=\vbox to
                 \ht0{}\,\box1^{#1}\!#2}

\newcommand{\lMod}[1]{\llap{\phantom{|}}_{#1}\catfont{Mod}}

\newcommand{\lindex}[2]{\llap{\phantom{|}}_{#1}{#2}}

\newcommand{\MackfR}{\mathrm{Mack}_R^f}
\newcommand{\MackfmuR}{\mathrm{Mack}_R^{f,\mu}}

\newcommand{\myiso}{\buildrel\sim\over\to}

\newcommand{\Ob}{\mathrm{Ob}}

\newcommand{\phiD}{\lexp{\varphi}{D}}

\newcommand{\qed}{\nobreak\hfill
                  \vbox{\hrule\hbox{\vrule\hbox to 5pt
                  {\vbox to 8pt{\vfil}\hfil}\vrule}\hrule}}
\newcommand{\Rab}{R^{\mathrm{ab}}}
\newcommand{\Res}{\mathrm{Res}}
\newcommand{\ResfR}{\mathrm{Res}_R^f}
\newcommand{\ResfmuR}{\mathrm{Res}_R^{f,\mu}}
\newcommand{\res}{\mathrm{res}}

\newcommand{\sigmatilde}{\tilde{\sigma}}
\newcommand{\Sigmatilde}{\tilde{\Sigma}}
\newcommand{\Sp}{\mathrm{Sp}}
\newcommand{\Sptilde}{\widetilde{\Sp}}

\newcommand{\tautilde}{\tilde{\tau}}

\newcommand{\xL}{\lexp{x}{L}}

\newcommand{\ZZ}{\mathbb{Z}}

\title{The $-_+$ and $-^+$ Constructions for Biset Functors
\footnote{MR Subject Classification 19A22, 20C15, 20C20}
\footnote{
Partially supported by Projects 2015 CN-15-43 UC MEXUS-CONACYT Collaborative Research Grants "Representation rings of finite groups", PAPIIT IN101416 Anillo global de representaciones y funtores asociados, and by the CIC at the UMSNH}
}
\author{\small Robert Boltje\\
        \small Department of Mathematics\\
        \small University of California\\
        \small Santa Cruz, CA 95064\\
        \small U.S.A.\\
        \small boltje@math.ucsc.edu %
        \and
        \small Gerardo Raggi-C\'ardenas\\
        \small Centro de Ciencias Matem\'aticas\\
        \small UNAM\\
        \small Morelia, Mich\\
        \small MEXICO\\
        \small raggi@matmor.unam.mx
        \and
        \small Luis Valero-Elizondo\\
        \small Facultad de Cs Fis-Mat\\
        \small Universidad Michoacana\\
        \small Morelia, Mich\\
        \small MEXICO\\
        \small valero@fismat.umich.mx
}
\date{May 23, 2018}

\begin{document}

\sloppy

\maketitle


\begin{abstract}
In this article we define the $-_+$-construction and the $-^+$-construction, that was crucial in the theory of canonical induction formulas (see \cite{Boltje1998b}), in the setting of biset functors, thus providing the necessary framework to define and construct canonical induction formulas for representation rings that are most naturally viewed as biset functors. Additionally, this provides a unified approach to the study of a class of functors including the Burnside ring, the monomial Burnside ring and global representation ring.
\end{abstract}


\section{Introduction}

This paper can be considered as the first step to extend the framework of canonical induction formulas, introduced by the first author in \cite{Boltje1998b}, 
to the setting of biset functors, a notion introduced by Bouc, see \cite{Bouc2010a}. The most basic example of a canonical induction formula is the one for 
the character ring $R(H)$ of a finite group, which one can regard as a canonical section of $a_H\colon R(H)\to \Rab_+(H)$ of the 
natural map $b_H\colon \Rab_+(H)\to R(H)$. Here, $\Rab_+(H)$ is the free abelian group on $H$-conjugacy classes $[K,\varphi]_H$ of 
pairs $(K,\varphi)$, where $K\le H$ and $\varphi$ is a linear character of $K$, i.e., a character of degree $1$, and $b_H([K,\varphi]_H)=\ind_K^H(\varphi)$. 
In \cite{Boltje1998b}, both $R$ and $\Rab_+$ were considered as Mackey functors, where $H$ runs through all subgroups of a fixed finite group $G$. 
The maps $b_H$, $H\le G$, commute with conjugations, restrictions, and inductions, while the maps $a_H$, $H\le G$, commute only with conjugations and 
restrictions. The groups $\Rab_+(H)$, $H\le G$, were constructed via the $-_+$-construction from the groups $\Rab(H)$, $H\le G$. Here $\Rab(H)$ 
denotes the subgroups of $R(H)$ generated generated by linear characters of $H$. These groups allow conjugation and restriction maps, 
but no induction maps. Canonical induction formulas exist for a variety of representation rings: The character ring, the Brauer character ring, 
the group of projective representations, the trivial source ring, and the linear source ring, see \cite{Boltje1998b}. There were two major 
unsatisfactory aspects to these constructions: The first is the restriction to subgroups of a fixed finite group $G$, the second one is the 
disregard for other natural operations, in particular {\em inflations}. Both these aspects were a consequence of using the framework of Mackey functors, 
the main functorial setup for representation rings available at that time. Bouc's theory of biset functors removes both of these restrictions and 
gives additional freedom of choosing certain of the operations of restriction, induction, inflation, deflation, conjugations, isomorphism, 
or even more general sets of operators.

\smallskip
The aim of this paper is to introduce the $-_+$ and $-^+$-constructions within the framework of biset functors. The $-^+$-construction is an auxiliary construction which is crucial in the proofs and the understanding of the $-_+$-construction. This should also be of interest independent of the theory of canonical induction formulas, since the $_+$ -construction yields various important biset functors (see Example~\ref{ex Burnside examples}: The Burnside functor is the $-_+$-construction of the constant biset functor with values $\ZZ$. For an abelian group $A$, the monomial (also called $A$-fibered) Burnside ring is the $-_+$-construction applied to the biset functor mapping a finite group $G$ to the free $\ZZ$-module with basis $\Hom(G,A)$, see \cite{Dress}, \cite{Barker}, or \cite{BoltjeCoskun} for instance. This was used in the known examples of canonical induction formulas for the representation rings mentioned above, where $A$ is a subgroup of the unit group of an appropriate field. More recently, the second and the third author introduced the notion of the {\em global representation ring}, see \cite{RaggiValero2015}, a combination of the Burnside ring and the character ring. Again, this construction turns out to be the $-_+$-construction applied to the character ring functor. Proving statements in general about the $-_+$-construction thus has applications for a variety of interesting examples, and unifies previous proofs for these examples. In the case of the first example, the Burnside ring $B(G)$ is the $-_+$-construction of the constant functor $\ZZ$ and the $-^+$-constructions of the constant functor $\ZZ$ yields its ghost ring $(\prod_{H\le G}\ZZ)^G$, where the exponent $G$ denotes taking $G$-fixed points with respect to permuting the components according to the conjugation action of $G$ on its subgroups. One of the main tools to study the Burnside ring is the {\em mark homomorphism} $B(G)\to(\prod_{H\le G}\ZZ)^G$ introduced by Burnside. This feature generalizes to our set-up.

\smallskip
We start our axiomatic setup with the choice of a family $\calG$ of finite groups and for every $G,H\in\calG$, a choice of a set $\calS(G,H)$ of subgroups of $G\times H$, satisfying axioms that lead to a category $\calD$, a subcategory of Bouc's biset category $\calC$. This general setup accommodates for instance the representation theory of the symmetric groups, where only symmetric groups on finite sets, their Young subgroups, and only restrictions and inductions between them are of interest. In Section~2 we recall basic definitions and facts on bisets and biset functors. Section~3 defines how to construct from $(\calG,\calS)$ associated pairs $(\calG,\calS_+)$ and $(\calG,\calS^+)$, which lead to subcategories $\calD_+$ and $\calD^+$ of $\calC$. Sections~4 and 5 describe the construction of the biset functor $F_+$ on $\calD_+$ and the biset functor $F^+$ on $\calD^+$ associated to a biset functor $F$ on $\calD$. In Section~6 we define the mark morphism $F_+\to F^+$ as a natural transformation. We prove that, under certain conditions on the base ring, the mark morphism is injective, or even bijective. Section~7 deals with the situation that $F$ has a multiplicative structure, more precisely, that $F$ is a Green biset functor. We show that then also $F_+$ and $F^+$ inherit Green biset functor structures and that the mark morphism is multiplicative. We also show how the species of $F_+(G)$ are determined by the species of $F(H)$, for $H\le G$ with $H\in\calG$. Finally, in Section~9 we prove adjunction properties of the functor $F\mapsto F_+$ similar to the properties proven in \cite{Boltje1998b}.


\newcommand{\gam}[1]{\Gamma_{#1}}
\newcommand{\subplus}[1]{#1_+}
\newcommand{\Gammathree}[3]{\Gamma_{#1,#2}(#3)}
\newcommand{\Gammatwo}[2]{\Gamma_{#1}(#2)}
\newcommand{\Gammaone}[1]{\Gamma_{#1}}
\newcommand{\bisetprod}[3]{#1\times_{#2}#3}
\newcommand{\Gammaplus}[2]{\Gamma_{+}(#1)(#2)}
\newcommand{\axid}{(id)\xspace}
\newcommand{\axconj}{(conj)\xspace}
\newcommand{\axcomp}{(comp)\xspace}
\newcommand{\axcompind}{(comp-ind)\xspace}
\newcommand{\axcompres}{(comp-res)\xspace}
\newcommand{\axexres}{(res)\xspace}
\newcommand{\axnodef}{(nodefl)\xspace}
\newcommand{\axexind}{(ind)\xspace}

\section{Bisets and Biset Functors}

Throughout this section, let $R$ denote a commutative ring (associative with $1$). We recall the notions of bisets and biset functors from \cite{Bouc2010a}.

\begin{nothing} {\em $(G,H)$-bisets and $B(G,H)$.}\quad 
For finite groups $G$ and $H$, a {\em $(G,H)$-biset} is a finite set $U$ equipped with a left $G$-action and a right $H$-action which commute: $g(uh)=(gu)h$ for all $g\in G$, $h\in H$, $u\in U$. The $(G,H)$-bisets form a category, whose morphisms are the functions that are $G$-equivariant and $H$-equivariant. Denote by $B(G,H)$ the Grothendieck group of $(G,H)$-bisets with respect to coproducts (disjoint unions). Identifying $(G,H)$-bisets with left $(G\times H)$-sets via $(g,h)u=guh^{-1}$ for $g\in G$, $h\in H$, and $u\in U$, the abelian group $B(G,H)$ is free with {\em standard basis elements} $[G\times H/D]$, where $D\le G\times H$ runs through a set of representatives of the conjugacy classes of subgroups of $G\times H$. Here, for a $(G,H)$-biset $U$, $[U]$ denotes the associated element in the Grothendieck group.

\smallskip
If also $K$ is a finite group, if $U$ is a $(G,H)$-biset, and $V$ is an $(H,K)$-biset, then $U\times V$ is an $H$-set via $h(u,v):=(uh^{-1},hv)$, for $u\in U$, $v\in V$ and $h\in H$. The $H$-orbit of $(u,v)$ is denoted by $[u,_H v]$ (or just $[u,v]$ if there is no risk of confusion) and the set of $H$-orbits is denoted by $U\times_H V$. The latter is naturally a $(G,K)$-biset via $g[u, v]k:=[gu, vk]$ for $(u,v)\in U\times V$ and $(g,k)\in G\times K$. This construction induces a bilinear map
\begin{equation}\label{eqn ten of bisets}
  -\cdotH-\colon B(G,H)\times B(H,K)\to B(G,K)\,.
\end{equation}
\end{nothing}

\begin{nothing} {\em The biset category\, $\calC$ and the biset functor category $\calF_{\calD,R}$.}\quad
Let $\calC$ denote the following category. Its objects are the finite groups, and for finite groups $G$ and $H$, one sets $\Hom_{\calC}(H,G):=B(G,H)$. The identity morphism of $G$ is $[G]$, where $G$ is viewed as $(G,G)$-biset by left and right multiplication. For finite groups $G,H,K$, the composition in $\calC$ is defined by the map in (\ref{eqn ten of bisets}). 

\smallskip
For any subcategory $\calD$ of $\calC$, a {\em biset functor on $\calD$ over $R$} is an additive functor $\calF\colon\calD\to \lMod{R}$. The biset functors on $\calD$ over $R$ form an abelian category $\calF_{\calD,R}$ whose morphisms are the natural transformations between biset functors.
\end{nothing}

\begin{nothing} {\em Subgroups of $G\times H$.}\quad
Let $G$ and $H$ be finite groups and let $D\le G\times H$. We write $p_1\colon G\times H\to G$ and $p_2\colon G\times H\to H$ for the natural projection maps. Moreover, we set 
\begin{equation*}
  k_1(D):=\{g\in G\mid (g,1)\in D\}\quad\text{and}\quad k_2(D):=\{h\in H\mid (1,h)\in D\}\,.
\end{equation*}
Note that the projection maps $p_i$ induce isomorphisms $D/(k_1(D)\times k_2(D))\myiso p_i(D)/k_i(D)$, for $i=1,2$. The resulting canonical isomorphism $\eta_D\colon p_2(D)/k_2(D)\myiso p_1(D)/k_1(D)$ is characterized by $\eta(hk_2(D))=gk_2(D)$ if and only if $(g,h)\in D$, for $g\in p_1(D)$ and $h\in p_2(D)$.
\end{nothing}

\begin{nothing}\label{noth * product} {\em The $*$-product.}\quad
For subgroups $D\le G\times H$ and $E\le H\times K$, one sets
\begin{equation*}
  D*E:=\{(g,k)\in G\times K\mid \exists h\in H \text{ with } (g,h)\in D \text{ and } (h,k)\in E\}\le G\times K\,.
\end{equation*}
Note that this product is associative. For a subgroup $H_1\le H$, one defines $D*H_1\le G$ in a similar way, by identifying $H$ with $H\times \{1\}$ and $G$ with $G\times \{1\}$.
Note also that if $U$ is a $(G,H)$-biset, $V$ is an $(H,K)$-biset, and $(u,v)\in U\times V$, then 
\begin{equation}\label{eqn stab *}
  (G\times K)_{[u,_H v]} = (G\times H)_u * (H\times K)_v\,.
\end{equation}
\end{nothing}

The following theorem gives an explicit formula for the product in (\ref{eqn ten of bisets}) in terms of standard basis elements, see \cite[Lemma~2.3.24]{Bouc2010a}.

\begin{theorem}\label{thm Mackey formula}(\cite[Lemma~2.3.24]{Bouc2010a})
Assume that $G$, $H$, $K$ are finite groups and $D\le G\times H$, $E\le H\times K$ are subgroups. Then,
\begin{equation*}
  \left[\frac{G\times H}{D}\right]\cdotH\left[\frac{H\times K}{E}\right] = \coprod_{t\in [p_2(D)\backslash H/p_1(E)]}
  \left[\frac{G\times K}{D*\lexp{(t,1)}{E}}\right] \ \in B(G,K)\,.
\end{equation*}
\end{theorem}

\begin{nothing}{\em Elementary biset operations.}
Assume that $G$ is a finite group. For a subgroup $H\le G$ one sets
\begin{equation*}
  \res^G_H:=[\lindex{H}{G}_G]=\left[\frac{H\times G}{\Delta(H)}\right]\in B(H,G)\quad\text{and}\quad
  \ind_H^G:=[\lindex{G}{G}_H]=\left[\frac{G\times H}{\Delta(H)}\right]\in B(G,H)\,,
\end{equation*}
where $G$ is viewed as a $(G,H)$-biset and as an $(H,G)$-biset with via left and right multiplication and $\Delta(H):=\{(h,h)\mid h\in H\}$. For a normal subgroup $N\trianglelefteq G$ one sets
\begin{gather*}
  \infl^G_{G/N}:=[\lindex{G}{(G/N)}_{G/N}]=\left[\frac{G\times G/N}{\{(g,gN)\mid g\in G\}}\right]\in B(G,G/N)\quad\text{and}\\
  \defl^G_{G/N}:=[\lindex{G/N}{(G/N)}_G]=\left[\frac{G/N\times G}{\{(gN,g)\mid g\in G\}}\right]\in B(G/N,G)\,,
\end{gather*}
where $G/N$ is viewed as a $(G,G/N)$-biset and as a $(G/N,G)$-biset via the natural epimorphism $G\to G/N$ and left and right multiplication. Finally, if $\alpha\colon G\myiso G'$ is a group isomorphism, we set $\isog_\alpha:=[G']\in B(G',G)$, where $G'$ is considered as $(G',G)$-biset via $g'xg:=g'x\alpha(g)$ for $g',x\in G'$ and $g\in G$. The five elements defined above are referred to as {\em restriction}, {\em induction}, {\em inflation}, {\em deflation}, and {\em isogation}. When $F$ is a biset functor then these elements induce maps between the respective evaluations of $F$. \end{nothing}

For arbitrary finite groups $G$ and $H$, one has a canonical decomposition of a standard basis element $\left[\frac{G\times H}{D}\right]$ of $B(G,H)$ into five elementary bisets:

\begin{theorem}\label{thm elementary decomposition}(\cite[Lemma~2.3.26]{Bouc2010a})
Let $G$ and $H$ be finite groups and let $D\le G\times H$. Then, in the category $\calC$, the morphism $\left[\frac{G\times H}{D}\right]$ can be written as the following composition of elementary biset operations:
\begin{equation*}
   \ind_{p_1(D)}^G\circ\infl_{p_1(D)/k_1(D)}^{p_1(D)}\circ\isog_{\eta_D}\circ\defl^{p_2(D)}_{p_2(D)/k_2(D)}\circ\res^{H}_{p_2(D)}\,.
\end{equation*}
\end{theorem}


\section{The $-^+$ and $-_+$ Constructions on Subcategories $\calD$ of $\calC$}\label{sec D_+}

For a finite group $G$ we denote by $\Sigma(G)$ the set of all subgroups of $G$.

\begin{nothing}\label{noth GS} {\em The data $(\calG,\calS)$.}\quad
(a) Throughout this section we consider a class $\calG$ of finite groups and a family $\calS=(\calS(G,H))_{G,H\in\calG}$, with $\calS(G,H)\subseteq\Sigma(G\times H)$ for $G,H\in\calG$. We will assume throughout, that $\calG$ and $\calS$ satisfy the following axioms:

\smallskip
(i) For all $G\in\calG$ one has $\Delta(G)\in\calS(G,G)$.

\smallskip
(ii) For all $G,H\in\calG$ the set $\calS(G,H)$ is closed under $G\times H$-conjugation.

\smallskip
(iii) For all $G,H,K\in\calG$ and all $D\in\calS(G,H)$ and $E\in\calS(H,K)$ one has $D*E\in\calS(G,K)$.

For $G\in\calG$ we will denote by $\Sigma_\calG(G)$ the set of all subgroups $H$ of $G$ with $H\in\calG$.

\smallskip
(b) In the sequel we will also consider additional properties of $(\calG,\calS)$ that we will require as necessary. They are

\smallskip
(iv) For all $G,H\in \calG$, all $D\in\calS(G,H)$ and all $K\in\Sigma_{\calG}(H)$, one has $D*K\in\calG$ and $D*\Delta(K)\in\calS(D*K,K)$. Note that $D*\Delta(K)=D\cap (G\times K)$ and that $p_1(D*\Delta(K))=D*K$, but in general $p_2(D*\Delta(K))=p_2(D)\cap K$ can be a proper subgroup of $K$.

\smallskip
(v) For all $G\in\calG$ and all $H\in\Sigma_{\calG}(G)$, one has $\Delta(H)\in\calS(G,H)$. 

\smallskip
(vi) For all $G\in\calG$ and all $H\in\Sigma_{\calG}(G)$, one has $\Delta(H)\in\calS(H,G)$.

%
\smallskip
(vii) For all $G,H\in\calG$ and all $D\in\calS(G,H)$ one has $p_2(D)\in\calG$, and for all $K\in\Sigma_{\calG}(p_2(D))$ one has $D*K\in\calG$ and $D*\Delta(K)\in\calS(D*K,K)$. Note that $D*\Delta(K)=D\cap (G\times K)$, $p_1(D*\Delta(K))=D*K$ and $p_2(D*\Delta(K))=K$. Note also that this condition is symmetric (cf.~Proposition~\ref{prop D^+}).

\smallskip
(c) Assume that $(\calG,\calS)$ satisfies the additional axiom (iv) or (vii). If $G\in\calG$, $H\in \Sigma_{\calG}(G)$, and $g\in G$, then also $\lexp{g}{H}\in\calG$ and $\lexp{(g,1)}{\Delta(H)}\in\calS(\lexp{g}{H},H)$. In fact, $\lexp{(g,1)}{\Delta(G)}\in\calS(G,G)$ by (i) and (ii), and by (iv) or (vii) we obtain $\lexp{g}{H}=\lexp{(g,1)}{\Delta(G)}*H\in\calG$ and $\lexp{(g,1)}{\Delta(H)}=\lexp{(g,1)}{\Delta(G)}*\Delta(H) \in\calS(\lexp{g}{H},H)$.
\end{nothing} 

\begin{nothing}\label{noth D(G,S)} {\em The category $\calD=\calC(\calG,\calS)$.}\quad
Given $(\calG,\calS)$ as in \ref{noth GS}(a), we define the subcategory $\calD=\calC(\calG,\calS)$ of the biset category $\calC$ by $\Ob(\calD)=\calG$ and, for $G,H\in\calG$, 
\begin{equation*}
  \Hom_{\calD}(H,G) := \big\langle\left[\frac{G\times H}{D}\right] \mid D\in\calS(G,H)\big\rangle_{\ZZ} \subseteq B(G,H)\,.
\end{equation*}
Axioms (i)---(iii) in \ref{noth GS}(a) and Theorem~\ref{thm Mackey formula} imply immediately that this is in fact a subcategory of $\calC$. Note that Axiom~\ref{noth GS}(v) (resp.~(vi)), if valid, ensures that the category $\calD$ contains all possible inductions (resp.~restrictions).
\end{nothing}

\begin{definition}\label{def D_+}
Let $(\calG,\calS)$ be as in \ref{noth GS}(a). For $G,H\in\calG$ we define $\calS_+(G,H)$ as the set of all subgroups $D\le G\times H$ such that 
\begin{equation*}
  \calS_+(G,H):=\{ D\le G\times H\mid p_1(D)\in\calG \text{ and }D\in\calS(p_1(D),H)\}\,.
\end{equation*}
The following proposition shows that if $(\calG,\calS)$ satisfies also Axiom~(iv) then $(\calG,\calS_+)$ satisfies again the axioms (i)--(iii) in \ref{noth GS}(a), so that we obtain a category $\calD_+:=\calC(\calG,\calS_+)$ by \ref{noth D(G,S)}.
\end{definition}

Note that each $D\in\calS_+(G,H)$ can be written as $\Delta(p_1(D))*D$ with $\Delta(p_1(D))\le G\times p_1(D)$ and $D\in\calS(p_1(D),H)$. Thus, $\left[\frac{G\times H}{D}\right] = \ind_{p_1(D)}^G \circ \left[\frac{p_1(D)\times H}{D}\right]$ with $\left[\frac{p_1(D)\times H}{D}\right]\in\Hom_{\calD}(H,p_1(D))$.

\begin{proposition}\label{prop D_+}
Assume that $(\calG,\calS)$ is as in \ref{noth GS}(a), satisfying Axioms~(i)--(iii) and additionally Axiom~(iv).

\smallskip
{\rm (a)} For $G, H\in\calG$ and $D\le G\times H$, one has $D\in\calS_+(G,H)$ if and only if, for all $K\in\Sigma_{\calG}(H)$, one has $D*K\in\calG$ and $D*\Delta(K)\in\calS(D*K,K)$. In particular, by Axiom~(iv), $\calS(G,H)\subseteq\calS_+(G,H)$ for all $G,H\in\calG$.

\smallskip
{\rm (b)} $(\calG,\calS_+)$ satisfies the axioms (i)--(v) in \ref{noth GS}.

\smallskip
{\rm (c)} One has $\calD\subseteq \calD_+$, with equality if and only if $(\calG,\calS)$ satisfies Axiom~(v). In particular, by Part~(b), $(\calD_+)_+=\calD_+$, and $\calD_+$ is the smallest $\ZZ-linear$ subcategory of $\calC$ containing $\calD$ and all inductions $\ind_{H}^G$ with $G\in\calG$ and $H\in\Sigma_{\calG}(G)$.

\smallskip
{\rm (d)} Let $G,H\in\calG$ and $D\le G\times H$ with $p_1(D)=G$. Then $D\in\calS(G,H)$ if and only if $D\in\calS_+(G,H)$. In particular, $\calD_+$ contains a given restriction, inflation, or deflation if and only if $\calD$ does.
\end{proposition}

\begin{proof}
(a) Let $G,H\in\calG$ and $D\le G\times H$. First assume $D\in\calS_+(G,H)$ and $K\in\Sigma_{\calG}(H)$. Then $p_1(D)\in\calG$ and $D\in\calS(p_1(D),H)$ by the definition of $\calS_+(G,H)$. Applying Axiom~(iv) to $D\in\calS(p_1(D),H)$ and $K$, we obtain $D*K\in\calG$ and $D*\Delta(K)\in\calS(D*K,K)$. For the converse, apply the condition to $K=H$ and note that $D*H=p_1(D)$ and $D*\Delta(H)=D$.

\smallskip
(b) Clearly, $(\calG,\calS_+)$ satisfies Axiom~(i). 

To see that it satisfies Axiom~(ii), let $D\in\calS_+(G,H)$ and $(a,b)\in G\times H$. We have $p_1(\lexp{(a,b)}{D})=\lexp{a}{p_1(D)}\in\calG$ and $\lexp{(a,1)}{\Delta(p_1(D))}\in \calS(\lexp{a}{p_1(D)},p_1(D))$ by \ref{noth GS}(c). Thus, $\lexp{(a,1)}{D}=\lexp{(a,1)}{\Delta(p_1(D))}*D\in\calS(\lexp{a}{p_1(D)},H)$ by Axiom~(iii), and therefore, by Axiom~(ii), $\lexp{(a,b)}{D}=\lexp{(1,b)}{(}\lexp{(a,1)}{D}) \in\calS(\lexp{a}{p_1(D)},H) = \calS(p_1(\lexp{(a,b)}{D}),H)$. This shows that $(\calG,\calS_+)$ satisfies Axiom~(ii). 

Next we show that $(\calG,\calS_+)$ satisfies Axiom~(iii). Let $G,H,K\in\calG$ and $D\in\calS_+(G,H)$ and $E\in\calS_+(H,K)$. We will use the characterization from Part~(a). Let $L\in\Sigma_{\calG}(K)$. Since $(D*E)*L= D*(E*L)$, we see that $E*L\in\Sigma_\calG(H)$ and then $(D*E)*L\in\Sigma_\calG(G)$. Next we show that $(D*E)*\Delta(L)\in \calS(D*E*L,L)$. We have $(D*E)*\Delta(L) = (D*\Delta(E*L))*(E*\Delta(L))$. By Part~(a) we have $E*\Delta(L)\in \calS(E*L,L)$ and also $D*\Delta(E*L)\in\calS(D*E*L,E*L)$. Now Axiom~(iii) for $(\calG,\calS)$ implies that $(D*E)*\Delta(L)\in\calS(D*E*L,L)$.

Axiom~(iv) for $(\calG,\calS_+)$ follows immediately from Part~(a) and Axiom~(iv) for $(\calG,\calS)$. 

Axiom~(v) clearly holds for $(\calG,\calS_+)$ by the definition of $\calS_+$, since $\Delta(H)\in\calS(H,H)$, for all $H\in\calG$.

\smallskip
(c) Since $\calS(G,H)\subseteq\calS_+(G,H)$ by Part~(a), we have $\calD\subseteq\calD_+$. If $\calD_+=\calD$, then $\calS=\calS_+$ and, by Part~(b), $(\calG,\calS_+)$ satisfies Axiom~(v), so also $(\calG,\calS)$ does. Conversely, assume $(\calG,\calS)$ satisfies Axiom~(v), and let $D\in\calS_+(G,H)$, for some $G,H\in\calG$. Then $p_1(D)\in\calG$ and $D\in\calS(p_1(D),H)$. Moreover $\Delta(p_1(D))\in\calS(G,p_1(D))$ by Axiom~(v) for $(\calG,\calS)$. By Axiom~(iii) for $(\calG,\calS)$, also $D=\Delta(p_1(D))*D\in \calS(G,H)$. The last statement in (c) follows immediately from the paragraph following Definition~\ref{def D_+}.

\smallskip
(d) This follows immediately from the definition of $\calS_+$.
\end{proof}

\begin{definition}\label{def D^+}
Let $(\calG,\calS)$ be as in \ref{noth GS}(a). For $G,H\in\calG$ we define 
\begin{equation*}
  \calS^+(G,H):=\{D\le G\times H\mid \text{$p_1(D)\in \calG$, $p_2(D)\in\calG$, and $D\in\calS(p_1(D),p_2(D))$}\}\,.
\end{equation*}
The following proposition shows that if $(\calG,\calS)$ also satisfies Axiom~(vii) then $(\calG,\calS^+)$ satisfies again Axioms (i)--(iii) in \ref{noth GS}(a), so that we obtain a category $\calD^+:=\calC(\calG,\calS^+)$ by \ref{noth D(G,S)}.
\end{definition}

Note that each $D\in\calS^+(G,H)$ can be written as $\Delta(p_1(D))*D*\Delta(p_2(D))$ with $\Delta(p_1(D))\le G\times p_1(D)$, $D\in\calS(p_1(D),p_2(D))$, and $\Delta(p_2(D))\in\calS(p_2(D),H)$. Thus, $\left[\frac{G\times H}{D}\right] = \ind_{p_1(D)}^G \circ \left[\frac{p_1(D)\times p_2(D)}{D}\right]\circ\res^H_{p_2(D)}$ with $\left[\frac{p_1(D)\times p_2(D)}{D}\right]\in\Hom_{\calD}(p_2(D),p_1(D))$.

\begin{proposition}\label{prop D^+}
Assume that $(\calG,\calS)$ is as in \ref{noth GS}(a), satisfying Axioms~(i)--(iii) and additionally Axiom~(vii).

\smallskip
{\rm (a)} For $G,H\in G$ and $D\le G\times H$, the following are equivalent:

\smallskip
{\rm (i)} $D\in\calS^+(G,H)$.

\smallskip 
{\rm (ii)} $p_2(D)\in\calG$ and for all $L\in\Sigma_\calG(p_2(D))$ one has $D*L\in\calG$ and $D*\Delta(L)\in\calS(L*D,L)$. 

\smallskip
{\rm (iii)} $p_1(D)\in\calG$ and for all $K\in\Sigma_\calG(p_1(D))$ one has $K*D\in\calG$ and $\Delta(K)*D\in\calS(K,K*D)$.

\smallskip
In particular, $\calS(G,H)\subseteq\calS^+(G,H)$ for all $G,H\in\calG$.

\smallskip
{\rm (b)} $(\calG,\calS^+)$ satisfies Axioms~(i)--(iii), (v), (vi), and (vii).

\smallskip
{\rm (c)} One has $\calD\subseteq\calD^+$, with equality if and only if $(\calG,\calS)$ satisfies Axioms~(v) and (vi). In particular, by Part~(b), $(\calD^+)^+=\calD^+$ and $\calD^+$ is the smallest $\ZZ$-linear subcategory of $\calC$ containing $\calD$ and all inductions $\ind_H^G$ and restrictions $\res^G_H$ with $G\in\calG$ and $H\in\Sigma_\calG(G)$.

\smallskip
{\rm (d)} Let $G,H\in\calG$ and $D\le G\times H$ with $p_1(D)=G$ and $p_2(D)=H$. Then $D\in\calS(G,H)$ if and only if $D\in\calS^+(G,H)$. In particular, $\calD^+$ contains a given inflation or deflation if and only if $\calD$ does.
\end{proposition}

\begin{proof}
(a) Let $D\le G\times H$. First assume that $D\in\calS^+(G,H)$. Then, by definition, $p_2(D)\in\calG$. Further assume that $L\in\Sigma_\calG(p_2(D))$. Then, by Axiom~(vii) for $(\calG,\calS)$ and $D$, we have $D*L\in\calG$ and $D*\Delta(L)\in\calS(D*L,L)$. Thus (ii) holds. Conversely, assume that $D$ satisfies the condition in (ii). Then $p_2(D)\in\calG$. Applying the condition in (ii) to $L:=p_2(D)$, we obtain $p_1(D)=D*p_2(D)\in\calG$ and $D=D*\Delta(p_2(D))\in\calS(p_1(D),p_2(D))$. Thus, (i) holds. Since the definition of $\calS^+$ is symmetric, the equivalence of (ii) and (iii) is proved in the same way. The last statement in (a) follows immediately from this characterization and from $(\calG,\calS)$ satisfying Axiom~(vii).

\smallskip
(b) Axioms (i)--(iii) and (v) follow from the same arguments as in the proof of \ref{prop D_+}. Axiom~(vi) follows immediately from the definition of $\calS^+$, and Axiom~(vii) for $(\calG,\calS^+)$ follows from Part~(a) and Axiom~(vii) for $(\calG,\calS)$.

\smallskip
(c) This is proved in a similar way as Part~(c) of Proposition~\ref{prop D_+} using the decomposition of $\left[\frac{G\times H}{D}\right]$ in the paragraph following Definition~\ref{def D^+}.

\smallskip
(d) This follows immediately from the definition of $\calS^+(G,H)$.
\end{proof}

We will need the following lemma in the subsequent example. Its proof is an easy verification.

\begin{lemma}\label{lem conditions and *}
Let $G,H,K$ be finite groups and let $D\le G\times H$ and $E\le H\times K$.

\smallskip
{\rm (a)} Let $i\in\{1,2\}$. If $k_i(D)=\{1\}$ and $k_i(E)=\{1\}$ then $k_i(D*E)=\{1\}$.

\smallskip
{\rm (b)} If $p_1(D)=G$ and $p_1(E)=H$ then $p_1(D*E)=G$. If $p_2(D)=H$ and $p_2(E)=K$ then $p_2(D*E)=K$.
\end{lemma}

\begin{example}\label{ex D_+ and D^+}
We say that a subgroup $D$ of $G\times H$ satisfies condition $k_1$ (resp.~$k_2$, resp.~$p_1$, resp.~$p_2$) if $k_1(D)=\{1\}$ (resp.~$k_2(D)=\{1\}$, resp.~$p_1(D)=G$, resp.~$p_2(D)=H$). For any subset $C$ of $\{k_1,k_2,p_1,p_2\}$, let $\calS_C(G,H)$ denote the set of subgroups $D$ of $G\times H$ which satisfy all the conditions in $C$.

Clearly, for any class $\calG$ of finite groups and any choice of $C$, $(\calG,\calS_C)$ satisfies Axioms (i) and (ii) in \ref{noth GS}(a). Moreover, by Lemma~\ref{lem conditions and *}, $(\calG,\calS_C)$ also satisfies Axiom~(iii) of \ref{noth GS}(a). Thus, by \ref{noth D(G,S)}, we obtain a subcategory $\calD_C(\calG):=\calD(\calG,\calS_C)$ of $\calC$.

Assume from now on that $\calG$ is the class of all finite groups and write $\calD_C$ for $\calD_C(\calG)$. It is easy to verify, using Theorems~\ref{thm Mackey formula} and \ref{thm elementary decomposition}, that the following holds: $\calD_C$ is the subcategory of $\calC$ generated by all isogations together with all the elementary operations contained in $E\subseteq\{\res,\ind,\infl,\defl\}$, where $E$ contains $\res$ (resp.~$\ind$, resp.~$\infl$, resp.~$\defl$) if and only if $C$ does {\em not} contain $p_2$ (resp.~$p_1$, resp.~$k_1$, resp.~$k_2$). It is also an easy verification that $(\calG,\calS_C)$ satisfies Axioms (iv) and (vii) of \ref{noth GS}(a). Therefore, by Proposition~\ref{prop D_+}(c), we obtain $(\calD_C)_+=\calD_{C\smallsetminus\{p_1\}}$, or equivalently that $(\calD_C)_+$ is the subcategory of $\calC$ generated by all the elementary operations that are already contained in $\calD$ together with all inductions.
Similarly, by Proposition~\ref{prop D^+}(c), we obtain $(\calD_C)^+=\calD_{C\smallsetminus\{p_1,p_2\}}$, or equivalently that $(\calD_C)^+$ is generated by all the elementary operations that are already contained in $\calD$ together with all inductions and all restrictions.
\end{example}

\begin{example}\label{ex D_+ and D^+ local}
Let $G$ be a finite group and let $\calG_G$ denote the set of all subgroups of $G$. Replacing all isogations with just conjugation isomorphisms induced by elements in $G$ and using a similar notation as in Example~\ref{ex D_+ and D^+}, we obtain three subcategories $\calD_{G,C}$ of $\calC$ for the three choices of subsets $C=\{k_1,k_2,p_1,p_2\}$, $C=\{k_1,k_2, p_1\}$, and $C=\{k_1,k_2\}$, whose objects are the elements of $\calG_G$ and whose morphisms are generated by all conjugations (resp.~all conjugations and restrictions, resp.~all conjugations, restrictions and inductions). For a commutative ring $R$, the functor categories $\calF_{\calD_{G,C},R}$ for these three choices are then very closely related to the conjugation functors (resp.~restriction functors, resp.~Mackey functors) on $G$ over $R$, as defined in \cite[Def.~1.1]{Boltje1998b}. In fact, the categories $\calF_{\calD_{G,C},R}$ are precisely the subcategories of the conjugation functors, restriction functors and Mackey functors in \cite{Boltje1998b}, satisfying the additional axiom that for $H\le G$ and $g\in C_G(H)$, the conjugation map induced by $g$ on the evaluation at $H$ is the identity. We will denote these categories by $\ConfR(G)$, $\ResfR(G)$, and $\MackfR(G)$, respectively. See \cite{Bouc2015} for a more detailed discussion on {\em fused} Mackey functors versus Mackey functors that elaborates on this aspect. We are grateful to Serge Bouc for pointing out this difference.
\end{example}


\section{The Functor $-_+\colon \calF_{\calD,R}\to\calF_{\calD_+,R}$}\label{sec F_+}

Throughout this section, let $R$ denote a commutative ring and let $(\calG,\calS)$ be as in \ref{noth GS}(a), satisfying Axioms~(i)--(iv). Let $\calD:=\calC(\calG,\calS)$ (cf.~\ref{noth D(G,S)}) and let $\calS_+$ and $\calD_+$ be defined as in Definition~\ref{def D_+}. The goal of this section is the construction of a functor $-_+\colon \calF_{\calD,R}\to\calF_{\calD_+,R}$ that generalizes the construction in \cite[Section~2]{Boltje1998b} in the situation of Example~\ref{ex D_+ and D^+ local}.

\begin{notation}\label{not exponential}
For a finite group $H$, an $H$-set $X$ and an element $x\in X$, we denote by $H_x:=\{h\in H\mid hx=x\}$ the stabilizer of $x$ in $H$. If also $G$ is a finite group and $U$ is a $(G,H)$-biset, then, for $u\in U$ and $K\le H$, we set
\begin{equation*}
  \lexp{u}{K}:=\{g\in G\mid \exists k\in K\colon gu=uk\}\,.
\end{equation*}
It is easy to verify that 
\begin{equation*}
  (G\times K)_u = (G\times H)_u*\Delta(K) = (G\times H)_u\cap(G\times K)\,,\quad   \lexp{u}{K}=p_1((G\times K)_u) = (G\times H)_u*K\,,
\end{equation*}
and that, for an $H$-set $X$ and $x\in X$, one has
\begin{equation*}
  G_{[u, x]}=\lexp{u}{H_x} = (G\times H)_u*H_x\,.
\end{equation*}
\end{notation}

\begin{nothing} {\em The category $\Gamma_F(G)$.}\quad 
Let $G\in\calG$, let $X$ be a finite $G$-set such that $G_x\in\calG$ for all $x\in X$, and let $F\in\calF_{\calD,R}$. A {\em section of $F$ over $X$} is a function $s\colon X\to \bigoplus_{x\in X} F(G_x)$ such that $s(x)\in F(G_x)$ for all $x\in X$. These sections form an $R$-module via point-wise constructions. The group $G$ acts $R$-linearly on the set of these sections by $(g\cdot s)(x):=\lexp{g}{s(g^{-1}x})$, where $\lexp{g}{m}:=F(\isog_{c_g})(m)$, for $g\in G$ and $m\in F(G_{g^{-1}x})$, and $c_g\colon G_{g^{-1}x}\to G_x$ is the conjugation isomorphism mapping $h\in G_{g^{-1}x}$ to $ghg^{-1}\in G_x$. A section $s$ of $G$ over $F$ is called {\em $G$-equivariant} if $g\cdot s=s$ for all $g\in G$.

\smallskip
For $G\in \calG$ and $F$ as above, we denote by $\Gamma_F(G)$ the category whose objects are the pairs $(X,s)$, where $X$ is a finite $G$-set  such that $G_x\in \calG$, for all $x\in X$, and $s$ is a $G$-equivariant section of $F$ over $X$. A morphism $\alpha\colon (X,s)\to (Y,t)$ in $\Gamma_F(G)$ is a morphism of $G$-sets $\alpha\colon X\to Y$ such that for all $x\in X$, one has $G_{\alpha(x)}=G_x$ and $t(\alpha(x))=s(x)$.
\end{nothing}

\begin{nothing} {\em The functor $\Gamma_F(U)\colon\Gamma_F(H)\to\Gamma_F(G)$.}\quad
Let $G,H\in\calG$, let $F\in\calF_{\calD,R}$, and let $U$ be a $(G,H)$-biset with $(G\times H)_u\in \calS_+(G,H)$ for all $u\in U$. We define a functor 
\begin{equation*}
\Gamma_F(U)\colon\Gamma_F(H)\to\Gamma_{F}(G)\,,\quad (X,s)\mapsto(U\times_H X,U(s))\,,
\end{equation*}
where
\begin{equation}\label{eqn U(s)}
  U(s)([u, x]):= F\left(\left[\frac{G_{[u,x]}\times H_x}{(G\times H_x)_u}\right]\right)(s(x))\,.
\end{equation}
For a morphism $\alpha\colon (X,s)\to(Y,t)$ in $\Gamma_F(H)$ we set $\Gamma_F(\alpha):=U\times_H\alpha$. For the rest of this subsection we show that these definitions are well-defined and yield a functor.

\smallskip
(a) First note that $H_x\in\calG$, for all $x\in X$, by the definition of $\Gamma_F(H)$ and that $G_{[u,x]}=(G\times H)_u * H_x\in\calG$, for all $u\in U$ and $x\in X$, by \ref{not exponential} and Proposition~\ref{prop D_+}(a). Moreover, $(G\times H_x)_u= (G\times H)_u*\Delta(H_x)\in\calS((G\times H)_u*H_x,H_x) = \calS(G_{[u,x]},H_x)$, again by \ref{not exponential} and Proposition~\ref{prop D_+}(a). Therefore, $F$ can be applied to the class of the biset in (\ref{eqn U(s)}).

\smallskip
(b) Next we show that the definition in (\ref{eqn U(s)}) does not change, when $(u,x)$ is replaced with $(uh^{-1},hx)$ for some $h\in H$. To see this, note that
\begin{equation}\label{eqn comp}
  F\left(\left[\frac{G_{[uh^{-1}, hx]}\times H_{hx}}{(G\times H_{hx})_{uh^{-1}}}\right]\right)(s(hx)) =
  F\left(\left[\frac{G_{[u, x]}\times H_{hx}}{(G\times H_x)_u*\Delta_{c_{h^{-1}}}(H_{hx})}\right]\right)(s(hx))\,,
\end{equation}
where $\Delta_{c_{h^{-1}}}(H_{hx}):=\{(h^{-1}h'h,h')\mid h'\in H_{hx}\}$, since $(G\times H_{hx})_{uh^{-1}} = 
(G\times H_x)_u*\Delta_{c_{h^{-1}}}(H_{hx})$, as a quick computation shows. Further,
\begin{equation*}
  \frac{G_{[u, x]}\times H_{hx}}{(G\times H_x)_u*\Delta_{c_{h^{-1}}}(H_{hx})} \cong 
  \frac{G_{[u, x]}\times H_x}{(G\times H_x)_u} \times_{H_x} \frac{H_x\times H_{hx}}{\Delta_{c_{h^{-1}}}(H_{hx})}\,,
\end{equation*}
so that, using the functor property of $F$, we may continue (\ref{eqn comp}) with
\begin{equation*}
  \cdots = \left(F\left(\left[\frac{G_{[u, x]}\times H_x}{(G\times H_x)_u}\right]\right)\circ 
  F(\isog_{c_{h^{-1}}})\right)(s(hx))
  = F\left(\left[\frac{G_{[u, x]}\times H_x}{(G\times H_x)_u}\right]\right)(\lexp{h^{-1}}{s(hx)})\,,
\end{equation*}
and arrive at the same expression, since $h\cdot s=s$.

\smallskip
(c) Next we show that the section $U(s)$ is $G$-equivariant. For $g\in G$ we have 
\begin{align*}
  (g\cdot U(s))[u, x] & = \lexp{g}{\bigl(U(s)([g^{-1}u,x])\bigr)} \\
  & = F\left(\left[\frac{G_{[u, x]} \times G_{[g^{-1}u, x]}}{\Delta_{c_g}(G_{[g^{-1}u, x]})}\times_{G_{[g^{-1}u, x]}}
          \frac{G_{[g^{-1}u, x]}\times H_x}{(G\times H_x)_{g^{-1}u}}\right]\right)(s(x)) \\
  & = F\left(\left[\frac{G_{[u, x]}\times H_x}{\Delta_{c_g}(G_{[g^{-1}u, x]})*(G\times H_x)_{g^{-1}u}}\right]\right)(s(x))\\
  & = F\left(\left[\frac{G_{[u, x]}\times H_x}{(G\times H_x)_u}\right]\right)(s(x))\,,
\end{align*}
since $\Delta_{c_g}(G_{[g^{-1}u, x]})*(G\times H_x)_{g^{-1}u} = (G\times H_x)_u$.

\smallskip
(d) Next let $(X,s)$ and $(Y,t)$ be objects in $\Gamma_F(H)$ and let $\alpha\colon X\to Y$ be a morphism of $H$-sets such that, for all $x\in X$, one has $H_{\alpha(x)}=H_x$ and $s(x)=t(\alpha(x))$. We need to show that $U\times_H \alpha$ is a morphism in $\Gamma_F(G)$ between $(U\times_H X,U(s))$ and $(U\times_H Y, U(t))$. Since $H_{\alpha(x)}=H_x$ we also have $G_{[u, x]}=G_{[u, \alpha(x)]}$ for all $u$ in $U$ and all $x\in X$ (see \ref{not exponential}). Moreover, for all $u\in U$ and $x\in X$, we have
\begin{align*}
  U(t)([u, \alpha(x)]) 
  & = F\left(\left[\frac{G_{[u,\alpha(x)]}\times H_{\alpha(x)}}{(G\times H_{\alpha(x)})_u}\right]\right)
  \bigl(t(\alpha(x))\bigr)\\
  & = F\left(\left[\frac{G_{[u,x]}\times H_x}{(G\times H_x)_u}\right]\right) \bigl(s(x)\bigr) = U(s)([u, x])\,,
\end{align*}
so that $U\times_H\alpha$ is a morphism in $\Gamma_F(G)$.

\smallskip
(e) Since $U\times_H -$ preserves compositions and the identity, we have verified that the above definitions yield a well-defined functor $\Gamma_F(U)\colon \Gamma_F(H)\to\Gamma_F(G)$.       
\end{nothing}

\begin{proposition}\label{prop functorial in U}
Let $G$, $H$, and $K$ be finite groups in $\calG$, let $U$ be a finite $(G,H)$-biset such that $(G\times H)_u\in\calS_+(G,H)$ for all $u\in U$, and let $V$ be a finite $(H,K)$-biset such that $(H\times K)_v\in\calS_+(H,K)$ for all $v\in V$. Then $(G\times K)_{[u,v]}\in \calS_+(G,K)$ for all $u\in U$ and $v\in V$. Moreover, the two functors $\Gamma_F(U\times _HV)$ and $\Gamma_F(U)\circ\Gamma_F(V)$ from $\Gamma_F(K)$ to $\Gamma_F(G)$ are naturally isomorphic.
\end{proposition}

\begin{proof}
The first statement follows immediately from Equation~(\ref{eqn stab *}) and Axiom (iii) in \ref{noth GS}(a), which holds for $(\calG, \calS_+)$ by Proposition~\ref{prop D_+}(b).

The two functors applied to an object $(X,s)$ of $\Gamma_F(K)$ yield the objects 
\begin{equation*}
  ((U\times_H V)\times_K X, (U\times_H V)(s)) \text{ and }(U\times_H(V\times_K X),U(V(s)))
\end{equation*}
of $\Gamma_F(G)$. Let 
\begin{equation*}
  \alpha\colon (U\times_H V)\times_K X\myiso U\times_H (V\times_K X)\,,\quad [[u,_H v],_Kx]\mapsto [u,_H[v,_K x]]\,,
\end{equation*}
for $(u,v,x)\in U\times V\times X$, be the usual natural isomorphism of $G$-sets. In order to show that $\alpha$ is a morphism between $\Gamma_F(U\times_H V)(X,s)$ and $(\Gamma_F(U)\circ\Gamma_F(V))(X,s)$, we still need to show that $G_{[[u,_H v],_Kx]}=G_{[u,_H[v,_K x]]}$ and that $((U\times_H V)(s))([[u,_H v],_Kx]) = U(V(s))([u,_H[v,_K x]])$, for all $(u,v,x)\in U\times V\times X$. The first statement is immediate, since $\alpha$ is an isomorphism of $G$-sets. For the second statement, we evaluate both sides. The left hand side equals
\begin{equation}\label{eqn evaluation}
  F\left(\left[\frac{G_{[[u,_H v],_K x]}\times K_x}{(G\times K_x)_{[u,_H v]}}\right]\right) (s(x))\,,
\end{equation}
and the right hand side equals
\begin{align*}
  &\quad F\left(\left[\frac{G_{[u,_H[v,_K x]]}\times H_{[v,_K x]}}{(G\times H_{[v,_K x]})_u}\right]\right)(V(s)(x)) \\
  & = \left(F\left(\left[\frac{G_{[u,_H[v,_K x]]}\times H_{[v,_K x]}}{(G\times H_{[v,_K x]})_u}\right]\right) \circ
       F\left(\left[\frac{H_{[v,_K x]}\times K_x}{(H\times K_x)_v}\right]\right)\right)(s(x)) \\
  & = F\left(\left[\frac{G_{[u,_H[v,_K x]]}\times H_{[v,_K x]}}{(G\times H_{[v,_K x]})_u}\times_{H_{[v,_K x]}}
           \frac{H_{[v,_K x]}\times K_x}{(H\times K_x)_v}\right]\right)(s(x))\,.
\end{align*}
A quick computation shows that $p_1((H\times K_x)_v) = H_{[v,_K x]}$. Thus, Theorem~\ref{thm Mackey formula} implies that the last expression is equal to
\begin{equation*}
  F\left(\left[\frac{G_{[u,_H[v,_K x]]}\times K_x}{(G\times H_{[v,_K x]})_u*(H\times K_x)_v}\right]\right)(s(x))\,,
\end{equation*}
which coincides with the element in (\ref{eqn evaluation}), since $G_{[u,_H[v,_K x]]}=G_{[[u,_H v],_Kx]}$ and $(G\times H_{[v,_K x]})_u*(H\times K_x)_v = (G\times K_x)_{[u,_H v]}$. In fact, the first equation was established before and the second equation is an easy verification.
\end{proof}

\begin{definition}\label{def 2 operations}
For $G\in\calG$ and $F\in\calF_{\calD,R}$, we define the following two operations in the category $\Gamma_F(G)$.

\smallskip
(a) Given two arbitrary objects $(X,s)$ and $(Y,t)$ in $\Gamma_F(G)$, their {\em coproduct} $(X,s)\coprod(Y,t)$ is defined as $(X\coprod Y, s\coprod t)$, where $X\coprod Y$ is the disjoint union of $X$ and $Y$ and $s\coprod t$ is the section on $X\coprod Y$ which is defined on $X$ as $s$ and on $Y$ as $t$. This construction together with the obvious inclusions from $X$ and $Y$ to $X\coprod Y$ is also a  categorial coproduct in $\Gamma_F(G)$.

\smallskip
(b) Given two objects $(X,s)$ and $(X,t)$ with the same underlying $G$-set $X$, we have an object $(X,s+t)$, where $s+t$ is the pointwise sum of the two sections $s$ and $t$.
\end{definition}

\begin{definition}\label{def F_+}
Let $G\in\calG$ and $F\in\calF_{\calD,R}$. We define the abelian group $F_+(G)$ as the free abelian group on the set of isomorphism classes $\{X,s\}$ of objects $(X,s)$ of $\Gamma_F(G)$ modulo the subgroup generated by all elements of the form
\begin{equation}\label{eqn relations}
  \{X\coprod Y, s\coprod t\}-\{X,s\}-\{Y,t\} \quad\text{and}\quad \{X,s+u\}-\{X,s\}-\{X,u\}\,,
\end{equation}
where $(X,s),(X,u),(Y,t)$ are objects of $\Gamma_F(G)$. We will denote the coset of $\{X,s\}$ in $F_+(G)$ by $[X,s]$.
\end{definition}

Note that the above free abelian group is an $R$-module, via $r\{X,s\}:=\{X,rs\}$, using the point-wise $R$-module structure of sections on a fixed $G$-set $X$, and note that the subgroup generated by the elements in (\ref{eqn relations}) is an $R$-submodule. Thus, $F_+(G)$ has a natural $R$-module structure.

Note also that by the first type of relations in (\ref{eqn relations}), the classes of the elements $\{G/H,s\}$, where $H\in\Sigma_{\calG}(G)$, and where $s$ is a $G$-equivariant section of $F$ over the $G$-set $G/H$, form a generating set of the abelian group $F_+(G)$. 

For every element $a\in F(H)$, there exists a unique $G$-equivariant section $s_a$ of $F$ over $G/H$ with $s(H)=a$. We abbreviate the class of $\{G/H,s_a\}$ by $[H,a]_G\in F_+(G)$.

\begin{theorem}\label{thm F_+}
Let $R$ be a commutative ring, let $(\calG,\calS)$ be a as in \ref{noth GS}(a) satisfying Axioms~(i)--(iv), set $\calD:=\calC(\calG,\calS)$ and $\calD_+:=\calC(\calG,\calS)$, and let $F\in\calF_{\calD,R}$ be a biset functor on $\calD$ over $R$.

\smallskip
{\rm (a)} Mapping a finite group $G$ to the $R$-module $F_+(G)$ and an element $[U]\in B(G,H)$, where $U$ is a finite $(G,H)$-biset with point stabilizers in $\calS_+(G,H)$, to the $R$-linear map $F_+([U])\colon F_+(H)\to F_+(G)$, induced by the functor $\Gamma_F(U)\colon \Gamma_F(H)\to\Gamma_F(G)$, yields a biset functor $F_+\in\calF_{\calD_+,R}$.

\smallskip
{\rm (b)} The association $F\mapsto F_+$ defines an $R$-linear functor $-_+\colon \calF_{\calD,R}\to\calF_{\calD_+,R}$.

\smallskip
{\rm (c)} For each $G\in\calG$, one has an $R$-module isomorphism
\begin{equation*}
  F_+(G)\cong \left(\bigoplus_{H\in\Sigma_{\calG}(G)} F(H)\right)_G,
\end{equation*}
where the above direct sum $M:=\bigoplus_{H\in\Sigma_{\calG}(G)} F(H)$ is an $RG$-module with $g\in G$ mapping $a\in F(H)$ to $\lexp{g}{a}\in F(\lexp{g}{H})$, and where $M_G$ denotes the $G$-cofixed points of $M$, i.e., the $R$-module $M/IM$, where $I$ is the augmentation ideal of $RG$. The isomorphism associates the class of the element $a\in F(H)$, for $H\in\Sigma_\calG(G)$, to $[H,a]_G\in F_+(G)$.
\end{theorem}

\begin{proof}
(a) It is straightforward to show that if $U$ and $V$ are isomorphic $(G,H)$-bisets with stabilizers in $\calS_+(G,H)$, then the functors $\Gamma_F(U)$ and $\Gamma_F(V)$ are naturally isomorphic. Thus, every $(G,H)$-biset $U$ induces a well-defined group homomorphism $f_{[U]}$ between the free abelian groups associated to $\Gamma_F(H)$ and $\Gamma_F(G)$. Moreover, $f_{[U]}$ is an $R$-module homomorphism and maps elements of the type in (\ref{eqn relations}) to elements of the same type. Thus $f_{[U]}$ induces an $R$-module homomorphism $F_+([U])\colon F_+(H)\to F_+(G)$. It is also easy to see that $f_{[U\coprod V]}= f_{[U]}+f_{[V]}$. Thus, there is a unique group homomorphism $F_+\colon B(G,H)\to \Hom_R(F_+(H),F_+(G))$ with $F_+([U])$ as defined above for every $(G,H)$-biset $U$. Moreover, Proposition~\ref{prop functorial in U} implies that $F_+$ is in fact a functor in $\calF_{\calD_+,R}$.

\smallskip
(b) If $\varphi\colon F\to F'$ is a morphism between objects $F,F'\in\calF_{\calD,R}$, then, for every $G\in\calG$, one obtains an induced functor $\Gamma_F(G)\to\Gamma_{F'}(G)$, which maps an object $(X,s)$ to the object $(X,\phi(s))$, where $(\phi(s))(x):=\phi_{G_x}(s(x))$. This functor induces an $R$-linear map from $F_+(G)$ to $F'_+(G)$. Moreover, it is straightforward to verify that this construction respects compositions, sends the identity morphism to the identity morphism, and is $R$-linear.

\smallskip
(c) For every object $(X,s)$ of $\Gamma_F(G)$, we can define an element $\phi_G(X,s)\in (\bigoplus_{H\in\Sigma_\calG(G)} F(H))_G$ as the class of $\sum_{x\in [G\backslash X]} s(x)$, with $s(x)\in F(G_x)$, where $[G\backslash X]$ denotes a set of representatives of the $G$-orbits of $X$. Isomorphic objects lead to the same elements and the elements in (\ref{eqn relations}) are mapped to $0$. Altogether we obtain an $R$-linear map $\phi_G\colon F_+(G)\to (\bigoplus_{H\in\Sigma_\calG(G)} F(H))_G$, which maps $[H,a]_G$ to the class of the element $a\in F(H)$. Conversely, for every $H\in\Sigma_\calG(G)$, one has a map $F(H)\to F_+(G)$, $a\mapsto[H,a]_G$. Since $[\lexp{g}{H},\lexp{g}{a}]_G=[H,a]_G$, for all $g\in G$, one obtains a well-defined map in the opposite direction which maps the class of $a\in F(H)$ to $[H,a]_G$. Thus, the two constructed maps are inverse to each other.
\end{proof}

\begin{remark}
One could have defined $F_+(G)$ directly as $(\bigoplus_{H\in\Sigma_{\calG}(G)} F(H))_G$. But then a proof that $F\mapsto F_+$ yields a functor from $\calF_{\calD,R}$ to $\calF_{\calD_+,R}$ would have been longer, even more technical, and uglier. The reason to introduce the category $\Gamma_F(G)$ was to give a more conceptual and also shorter proof.

Note that $\bigoplus_{H\in \Sigma_\calG(G)} F(H) = F(G)\oplus (\bigoplus_{G\neq H\in\Sigma_\calG(G)} F(H))$ is a decomposition into $RG$-submodules and that $G$ acts trivially on $F(G)$. Thus, with the identification in Theorem~\ref{thm F_+}(c), we obtain a decomposition
\begin{equation}\label{eqn F_+ decomp}
  F_+(G) = F(G) \oplus F_+^<(G)\,,
\end{equation}
where $F_+^<(G):=(\bigoplus_{G\neq H\in\Sigma_\calG(G)} F(H))_G$. We denote the corresponding projection by
\begin{equation}\label{eqn pi}
  \pi_G\colon F_+(G)\to F(G)\,.
\end{equation}
Note that for $H\in\Sigma_{\calG}(G)$ and $a\in F(H)$ one has $\pi_G([H,a]_G)=0$ if $H\neq G$ and $\pi_G([H,a]_G)=a$ if $H=G$.

For all practical purposes it is easier to use the version $F_+(G)=(\bigoplus_{H\in\Sigma_{\calG}(G)} F(H))_G$ and its generating elements $[H,a]_G$. Using the isomorphism between the two versions we can translate the general biset operation as follows: Assume that $G,H\in\calG$, $U$ is a $(G,H)$-biset with point stabilizers in $\calS_+(G,H)$, that $K\in\Sigma_{\calG}(H)$ and $a\in F(K)$. Then
\begin{equation*}
  F_+([U])([K,a]_H) = \sum_{[u, hK]\in[G\backslash(U\times_H (H/K))]}\bigl[ G_{[u,hK]},
      F\left(\left[\frac{G_{[u, hK]}\times \lexp{h}{K}}{(G\times \lexp{h}{K})_u}\right]\right) (\lexp{h}{a})\bigr]_G\,,
\end{equation*}
where $[G\backslash(U\times_H (H/K))]$ denotes a set of representatives of the $G$-orbits of $U\times_H (H/K)$.
In particular, if $D\in\calS_+(G,H)$, then
\begin{equation}\label{eqn explicit F_+}
  F_+\bigl(\bigl[\frac{G\times H}{D}\bigr]\bigr)([K,a]_H) = 
  \sum_{h\in [p_2(D)\backslash H/K]} \bigl[D*\lexp{h}{K}, 
  F\bigl(\bigl[\frac{D*\lexp{h}{K}\times \lexp{h}{K}}{D*\Delta(\lexp{h}{K})}\bigr]\bigr) (\lexp{h}{a}) \bigr]_G \,,
\end{equation}
where $[p_2(D)\backslash H/K]$ denotes a set of representatives of the double cosets of $G$ with respect to $p_2(D)$ and $K$. This formula specializes in the case of elementary bisets as follows:
\begin{itemize}
\item Assume that $G,G'\in\calG$, $\alpha\colon G\myiso G'$ is an isomorphism with $\Delta_\alpha(G)\in\calS_+(G',G)$, $K\in\Sigma_\calG(G)$ and $a\in F(K)$. Then
\begin{equation*}
  F_+(\isog_\alpha)([K,a]_G) = [\alpha(K),F(\isog_\alpha)(a)]_{G'}\,.
\end{equation*}
\item Assume that $G\in\calG$, $H\in\Sigma_{\calG}(G)$ with $\Delta(H)\in\calS_+(H,G)$, $K\in\Sigma_{\calG}(G)$, and that $a\in F(K)$. Then
\begin{equation*}
  F_+(\res^G_H)([K,a]_G) = \sum_{g\in[H\backslash G/K]} [H\cap\lexp{g}{K},F(\res^{\lexp{g}{K}}_{H\cap\lexp{g}{K}})(\lexp{g}{a})]_H\,.
\end{equation*}
\item Assume that $G\in\calG$, $H\in\Sigma_{\calG}(G)$ with $\Delta(H)\in\calS_+(G,H)$, $K\in\Sigma_{\calG}(H)$, and that $a\in F(K)$. Then
\begin{equation*}
  F_+(\ind_H^G)([K,a]_H) = [K,a]_G\,.
\end{equation*}
Note that even if $\calD$ contains inductions, the induction operations for $F_+$ are independent of the induction operations for $F$.
\item Assume that $G\in\calG$, $N\trianglelefteq G$ with $G/N\in\calG$ and $\{(g,gN)\mid g\in G\}\in\calS_+(G,G/N)$, $N\trianglelefteq K\le G$ with $K/N\in\calG$, and $a\in F(K/N)$. Then
\begin{equation*}
  F_+(\infl_{G/N}^G)([K/N,a]_{G/N}) = [K,F(\infl_{K/N}^K)(a)]_G\,.
\end{equation*}
\item Assume that $G\in\calG$, $N\trianglelefteq G$ with $G/N\in\calG$ and $\{(gN,g)\mid g\in G\}\in\calS_+(G/N,G)$, that $K\in\Sigma_\calG(G)$ and $a\in F(K)$. Then
\begin{equation*}
  F_+(\defl^G_{G/N})([K,a]_G) =  \left[KN/N,F\left(\left[\frac{(KN/N)\times K}{\{(kN,k)\mid k\in K\}}\right]\right)(a)\right]_{G/N}\,.
\end{equation*}
If also $K/(K\cap N)\in\calG$ and $\isog_\alpha$ and $\defl^K_{K/(K\cap N)}$ are in $\calD$, where $\alpha\colon K/K\cap N\myiso KN/N$ is the canonical isomorphism induced by the inclusion $K\le KN$, then we can write the above as
\begin{equation*}
   \bigl[KN/N, \bigl(F(\isog_\alpha)F(\defl^K_{K/(K\cap N)})\bigr)(a)\bigr]_{G/N}\,.
\end{equation*}
\end{itemize}
\end{remark}

\begin{example}\label{ex Burnside examples}
Assume that $\calG$ is the class of all finite groups and consider $\calS(G,H):=\calS_{\{p_1\}}(G,H)=(\{D\le G\times H\mid p_1(D)=G\}$, for $G,H\in\calG$. Then $\calD:=\calC(\calG,\calS) = \calD_{\{p_1\}}$, cf.~Example~\ref{ex D_+ and D^+}. Then $\calS_+(G,H)$ consists of all subgroups of $G\times H$ and $\calD_+=\calC$.

\smallskip
(a) Consider the constant functor $F\in\calF_{\calD,\ZZ}$ associating to each $G\in\calG$ the abelian group $\ZZ$ and to each $[(G\times H)/D]$ (with $G,H\in\calG$ and $D\le G\times H$ with $p_1(D)=G$) the identity map from $F(H)$ to $F(G)$. This is indeed a functor by Theorem~\ref{thm Mackey formula}.  Then $F_+$ identifies to the Burnside ring functor $B$ by mapping $[H,a]_G\in F_+(G)$ to $a[G/H]\in B(G)$, for $H\le G$ and $a\in\ZZ$.

\smallskip
(b) Let $A$ be an abelian group and define $F(G)$ as the free $\ZZ$-module with basis $\Hom(G,A)$. Then $F\in\calF_{\calD,\ZZ}$, if endowed with the obvious isogation, restriction, inflation and deflation operations. The resulting biset functor $F_+$ is isomorphic to the $A$-fibered Burnside ring functor $B^A$, see \cite{Dress}, \cite{Barker}, or \cite{BoltjeCoskun}. In particular, if $A=\CC^\times$ is the unit group of the field of complex numbers, we obtain $F(G)=\Rab(G)$ and $F_+(G)=\Rab_+(G)$ for any finite group $G$, cf.~\cite{Boltje1998b}.

\smallskip
(c) If $F=R$ denotes the character ring biset functor on $\calD$ over $\ZZ$, then $F_+$ is the global representation ring functor introduced in \cite{RaggiValero2015}.
\end{example}


\section{The Functor $-^+\colon \calF_{\calD,R}\to\calF_{\calD^+,R}$}\label{sec F^+}

Throughout this section, let $R$ be a commutative ring and let $(\calG,\calS)$ be as in \ref{noth GS}(a) satisfying the axioms (i)-(iii) and (vii) in \ref{noth GS}. Let $\calD:=\calC(\calG,\calS)$ (cf.~\ref{noth D(G,S)}) and let $\calS^+$ and $\calD^+$ be defined as in Definition~\ref{def D^+}. We additionally assume that for all $G,H\in\calG$ and $D\in\calS^+(G,H)$ one has $k_2(D)=\{1\}$, i.e., that $(G\times H)/D$ is free as $H$-set. This is equivalent to assuming the same for all $D\in\calS(G,H)$. In fact, $\calS\subseteq\calS^+$, by Proposition~\ref{prop D^+}(a), and if $D\in\calS^+(G,H)$, then $D\in\calS(p_1(D),p_2(D))$, by definition. Thus, $\calD$ and $\calD^+$ don't have deflations.

\smallskip
The goal of this section is to construct a functor $-^+\colon\calF_{\calD,R}\to\calF_{\calD^+,R}$ that generalizes the construction in \cite[Section~2]{Boltje1998b} in the situation of Example~\ref{ex D_+ and D^+ local}, up to the modification of working with {\em fused} Mackey functors and {\em fused} conjugation functors, cf.~\cite{Bouc2015}.

\begin{notation}\label{not K^u}
Similar to the notation in \ref{not exponential}, one defines for a finite $(G,H)$-biset $U$, $u\in U$, and $K\le G$,
\begin{equation*}
  K^u:=\{h\in G\mid \exists k\in K\colon (k,h)\in (K\times H)_u\}\le H\,.
\end{equation*}
Note that 
\begin{equation*}
  (K\times H)_u = \Delta(K)*(G\times H)_u\quad\text{and}\quad K^u=K*(G\times H)_u\,.
\end{equation*}
\end{notation}

\begin{definition}\label{def F^+}
For $G\in\calG$ we define
\begin{equation*}
  F^+(G):=\left(\prod_{H\in\Sigma_{\calG}(G)} F(H)\right)^G\,,
\end{equation*}
where the exponent $G$ denotes taking $G$-fixed points. That is, a tuple of elements $a_H\in F(H)$, for $H\in\Sigma_\calG(G)$, belongs to $F^+(G)$ if and only if $\lexp{g}{a_H}=a_{\lexp{g}{H}}$ for all $H\in\Sigma_{\calG}(G)$ and all $g\in G$. Note that if $\Sigmatilde_\calG(G)\subseteq\Sigma_\calG(G)$ is a set of representatives of the conjugacy classes of $\Sigma_\calG(H)$ then projection onto the components indexed by $\Sigmatilde_{\calG}(G)$ defines an $R$-module isomorphism 
\begin{equation}\label{eqn F_+ projection iso}
  F^+(G)\myiso\prod_{H\in\Sigmatilde_\calG(G)} F(H)^{N_G(H)}\,.
\end{equation}

\smallskip
For any $G,H\in\calG$ and any $(G,H)$-biset $U$ such that $(G\times H)_u\in\calS^+(G,H)$ for all $u\in U$, we define
\begin{align}\label{eqn F^+(U)}
  \notag F^+([U])=F^+(U)\colon F^+(H) & \to F^+(G)\,,\\
  (a_L)_{L\in\Sigma_{\calG}(H)} & \mapsto \left(\sum_{\substack{u\in[U/H]\\ K\le p_1((G\times H)_u)}}
           F\left(\left[\frac{K\times K^u}{(K\times H)_u}\right]\right)(a_{K^u})\right)_{K\in\Sigma_{\calG}(G)}\,,
\end{align}
where $[U/H]$ denotes a set of representatives of the $H$-orbits of $U$. The above expression is well-defined, since $K\in\calG$ implies $K^u=K*(G\times H)_u\in\calG$ and $(K\times H)_u= \Delta(K)*(G\times H)_u\in\calS(K,K^u)$, by \ref{not K^u} and Proposition~\ref{prop D^+}(a) (i)$\Rightarrow$(iii). Thus, $F$ can be applied to the class of the $(K,K^u)$-biset $(K\times K^u)/(K\times H)_u$. Finally, the expression on the right hand side in (\ref{eqn F^+(U)}) does not depend on the choice of representatives $[U/H]$ and neither on the choice of $U$ in its isomorphism class. Clearly, $F^+([U])$ is also an $R$-module homomorphism. Since the sum in (\ref{eqn F^+(U)}) is additive in $[U]$, we also obtain an induced group homomorphism
\begin{equation*}
  F^+\colon B(G,H)\to \Hom_R(F^+(H),F^+(G))\,.
\end{equation*}
\end{definition}

\begin{theorem}\label{thm F^+ is functor}
Let $(\calG,\calS)$ be as introduced at the beginning of the section, satisfying Axioms~(i)--(iii) and (vii) of \ref{noth GS} and assume that $k_2(D)=1$ for every $G,H\in\cal G$ and $D\in\calS(G,H)$. Set $\calD:=\calC(\calG,\calS)$ and $\calD^+=\calC(\calG,\calS^+)$. Then the constructions in Definition~\ref{def F^+} define an $R$-linear functor $-^+\colon\calF_{\calD,R}\to\calF_{\calD^+,R}$.
\end{theorem}

\begin{proof}
We first show that for all finite groups $G,H,K\in\calG$, all finite $(G,H)$-bisets $U$ with stabilizers in $\calS^+(G,H)$ and all finite $(H,K)$-bisets $V$ with stabilizers in $\calS^+(H,K)$, one has
\begin{equation}\label{eqn F^+ functor}
  F^+(U\times_H V) = F^+(U)\circ F^+(V)\colon F^+(K)\to F^+(G)\,.
\end{equation}
Let $c=(c_N)_{N\in\Sigma_\calG(K)}\in F^+(K)$. Then for every $L\in\Sigma_\calG(G)$, the $L$-component of $F^+(U\times_H V)(c)$  equals
\begin{equation}\label{eqn lhs}
  \sum_{\substack{[u,v]\in [(U\times_H V)/K]\\ L\in p_1((G\times K)_{[u, v]})}}
  F\left(\left[\frac{L\times L^{[u,v]}}{(L\times K)_{[u, v]}}\right]\right) (c_{L^{[u,v]}}) \,.
\end{equation}
Next we compute the right hand side of (\ref{eqn F^+ functor}) evaluated at $c$. Setting $b=(b_M)_{M\in\Sigma_\calG(H)}:=F^+(V)(c)$ and $a=(a_L)_{L\in\Sigma_{\calG}(G)}:= F^+(U)(b)$ we have, for every $M\in\Sigma_\calG(H)$,
\begin{equation*}
  b_M= \sum_{\substack{v\in [V/K]\\ M\le p_1((H\times K)_v)}} F\left(\left[\frac{M\times M^v}{(M\times K)_v}\right]\right)(c_{M^v})\,,
\end{equation*}
and for every $L\in\Sigma_\calG(G)$,
\begin{align*}
  a_L & = F^+(U)(b)= \sum_{\substack{u\in [U/H]\\ L\le p_1((G\times H)_u)}}
  F\left(\left[\frac{L\times L^u}{(L\times H)_u}\right]\right)(b_{L^u})\\\
  & = \sum_{\substack{u\in [U/H] \\ L\le p_1((G\times H)_u)}} \sum_{\substack{v\in [V/K] \\ L^u\le p_1((H\times K)_v)}}
  F\left(\left[\frac{L\times L^u}{(L\times H)_u}\times_{L^u}\frac{L^u\times (L^u)^v}{(L^u\times K)_v}\right]\right) (c_{(L^u)^v})\,.
\end{align*}
Since $U$ is right-free, for any sets of representatives $[U/H]$ and $[V/K]$, the set $\{[u,v]\mid (u,v)\in [U/H]\times [V/K]\}$ is a set of representatives of $(U\times_H V)/K$. Moreover,  by Lemma~\ref{lem double condition}(ii), the sum in (\ref{eqn lhs}) and the double sum in the last formula for $a_L$ run over the the same indexing set. Finally, we have $(L^u)^v=L^{[u,v]}$ by Lemma~\ref{lem double condition}(i) and $(L\times H)_u*(L^u\times K)_v = (L\times K)_{[u,v]}$, see Equation~(\ref{eqn stab *}). This finishes the proof of the equality in (\ref{eqn F^+ functor}). That $F^+$ preserves identity morphisms and is $R$-linear follows immediately from the definitions.
\end{proof}

The proof of the following lemma is an easy exercise left to the reader.

\begin{lemma}\label{lem double condition}
Let $U$ be a finite $(G,H)$-biset, let $V$ be a finite $(H,K)$-biset, let $(u,v)\in U\times V$ and let $L\le G$.

\smallskip
{\rm (a)} One has $(L^u)^v=L^{[u,_Hv]}$.

\smallskip
{\rm (b)} Assume that $U$ is right-free. Then $L\le p_1((G\times K)_{[u,v]})$ if and only if $L\le p_1((G\times H)_u)$ and $L^u\le p_1((H\times K)_v)$.
\end{lemma}


\section{The Mark Morphism}\label{sec mark morphism}

Let $(\calG,\calS)$ be as in \ref{noth GS} and assume that $(\calG,\calS)$ satisfies axioms (i)--(iv), (vi), and (vii) in \ref{noth GS} and the
condition $k_2$ in Example~\ref{ex D_+ and D^+}. Set $\calD:=\calC(\calG,\calS)$, $\calD_+:=\calC(\calG,\calS_+)$ and $\calD^+:=\calC(\calG,\calS^+)$. Then Propositions~\ref{prop D_+}(c) and \ref{prop D^+}(c) imply that $\calD_+=\calD^+$. We denote this category by $\calE$.
Let $F\in \calF_{\calD,R}$ and $G\in \calG$. We define an $R$-linear map
$m_{F,G}:=(\pi_H\circ F_+(\res^G_H))_{H\in\Sigma_\calG(G)} : F_+(G)\longrightarrow F^+(G)$, using the projection map $\pi_H\colon F_+(H)\to F(H)$ from (\ref{eqn pi}). Thus, $m_{F,G}$ maps the class $[X,s]\in F_+(G)$ of an object $(X,s)$ in $\Gamma_F(G)$ to the element $(a_L)_{L\in\Sigma_{\calG}(G)}\in F^+(G)$ with
$$a_L := \sum_{\substack{x\in X \\ L\leq G_x}}F(\res^{G_x}_L)(s(x))\,.$$
Note that for the definition of $F^+(G)$ and the map $m_{F,G}$ one does not need axiom~(vii) in \ref{noth GS} and the condition $k_2$ in Example~\ref{ex D_+ and D^+}. But it is necessary in the next theorem. In general, the map $m_{F,G}$ does not commute with deflations.


\begin{theorem}\label{thm mark diagram}
Assume that $(\calG,\calS)$ satisfies axioms (i)--(iv), (vi), and (vii) in \ref{noth GS} and the
condition $k_2$ in Example~\ref{ex D_+ and D^+}.
Let $G,H\in\calG$ and let $U$ be a finite $(G,H)$-biset with stabilizers in $\calS_+(G,H)=\calS^+(G,H)$. Then, for any $F\in\calF_{\calD,R}$, the diagram
$$\begin{CD}
F_+(H) @>m_{F,H}>> F^+(H)\\
@VF_+([U])VV @VVF^+([U])V\\
F_+(G) @>m_{F,G}>> F^+(G)
\end{CD}$$
commutes. In particular, the $R$-linear maps $m_{F,G}$, $G\in\calG$, define a morphism $m_F\colon F_+\to F^+$ in $\calF_{\calE,R}$ which we will call the {\em mark morphism} associated with $F$. 
\end{theorem}

%
%

\begin{proof}
Let $(X,s)$ be in $\Gamma_F(H)$, let $[X,s]\in F_+(H)$ be the associated class, and let $L\in\Sigma_\calG(G)$. Then the $L$-component of $(F^+([U])(m_{F,H}([X,s])))$ is equal to
\begin{align*}
  & \sum_{\substack{u\in [U/H] \\L\leq p_1((G\times H)_u)}}
  F\left(\left[\frac{L\times L^u}{(L\times H)_u}\right]\right)
  \Bigl(\sum_{\substack{x\in X\\L^{u}\leq H_x}}
  F(\res^{H_x}_{L^u})(s(x))\Bigr)\\
  = &  \sum_{\substack{u\in [U/H] \\x\in X\\L\leq p_1((G\times H)_u)\\L^{u}\leq H_x}}
  F\left(\left[\frac{L\times L^u}{(L\times H)_u}\times_{L^u}
  \frac{L^u\times H_x}{\Delta(L^u)}\right]\right)(s(x))\,.
\end{align*}
Moreover, the $L$-component of $m_{F,G}(F_+([U])([X,s]))=m_{F,G}([U\times_H X, U(s)])$ is equal to
\begin{align*}
& \sum_{\substack{[u,x]\in U\times_H X \\L\leq G_{[u,x]}}} F(\res^{G_{[u,x]}}_L)(U(s)([u,x]))\\
  = & \sum_{\substack{[u,x]\in U\times_H X \\L\leq G_{[u,x]}}}
F\left(\left[\frac{L\times G_{[u,x]}}{\Delta(L)}
\times_{G_{[u,x]}} \frac{G_{[u,x]}\times H_x}{(G\times H_x)_u}\right]
\right)(s(x))\\
  = & \sum_{\substack{[u,x]\in U\times_H X \\L\leq G_{[u,x]}}}
F\left(\left[ \frac{L\times H_x}{\Delta(L)* (G\times H_x)_u}\right]
\right)(s(x))\,.
\end{align*}
Since $U$ is right-free, the map $[U/H]\times X\to U\times_H X$, $(u,x)\mapsto [u,x]$, is bijective.
Moreover, Lemma~\ref{lem double condition}(b) applied to the case $K=\{1\}$ shows that the conditions for the
subgroup $L$ are equivalent in both sums. Thus, there is a bijection between the indexing sets of both sums. Finally,
notice that for all $L\in\Sigma_\calG(G)$, \ref{not exponential} and \ref{not K^u} imply that $\Delta(L)*(G\times H_x)_u = 
(L\times H_x)_u = (L\times H)_u * \Delta(L^u)$. This completes the proof.
\end{proof}

The following definition gives a map that is close to an inverse to the mark morphism.

\begin{definition}\label{def sigma}
Assume that $(\calG,\calS)$ satisfies axioms (i)--(iv), (vi) in \ref{noth GS} and set $\calD:=\calC(\calG,\calS)$.
Let $F\in\calF_{\calD,R}$ and $G\in\Sigma_\calG(G)$.
We define a map $n_{F,G}:F^+(G)\longrightarrow F_+(G)$ by
$$n_{F,G}((a_H)_{H\in \Sigma_{\calG}(G)}) = \sum_{\substack{L, K\in \Sigma_{\calG}(G)}}
|L| \mu(L,K) [L,\res^K_L(a_K)]\,.$$
Here, $\mu$ denotes the M\"obius function of the poset $\Sigma_\calG(G)$, ordered by inclusion.
\end{definition}

\begin{proposition}
Let $(\calG,\calS)$, $\calD$, $F$, and $G$ be as in Definition~\ref{def sigma}. Then $n_{F,G}\circ m_{F,G}=|G|\cdot\id_{F_+(G)}$ and $m_{F,G}\circ n_{F,G} = |G|\cdot\id_{F^+(G)}$.
\end{proposition}

\begin{proof}
Let $(X,s)\in\Gamma_F(G)$. We have
\begin{align*}
 (n_{F,G}\circ m_{F,G})([X,s])  
 = & \sum_{\substack{L, K\in\Sigma_{\calG}(G) }} | L | \mu(L,K) [L, F(\res^K_L)(m_{F,G}([X,s])_K)]_G\\
   = & \sum_{\substack{L,K\in\Sigma_\calG(G) }} 
  \sum_{\substack{x\in X \\ K\leq G_x }} | L | \mu(L,K) [L, F(\res^K_L \mathop{\cdot}\limits_K\res^{G_x}_K)(s(x))]_G\\
 = & \sum_{\substack{L\in\Sigma_\calG(G) }} 
  \sum_{\substack{x\in X}} | L |
  \sum_{\substack{K\in\Sigma_{\calG}(G)\\ L\leq K\leq G_x }}  \mu(L,K) [L, F(\res^{G_x}_L)(s(x))]_G\,,
\end{align*}
where $m_{F,G}([X,s])_K$ 
denotes the $K$-component of $m_{F,G}([X,s])$.
By the standard properties of $\mu$, the expression 
$$ \sum_{\substack{K\in\Sigma_\calG(G)\\ L\leq K\leq G_x }}  \mu(L,K)$$
equals 1 if $L=G_x$, and 0 otherwise. Therefore, the last expression equals

\begin{equation*}
  \sum_{\substack{L\leq G }} \sum_{\substack{x\in X \\ L = G_x}} | L | [L, s(x)]_G
  = \sum_{\substack{x\in X}} | G_x | [G_x, s(x)]_G 
  = \sum_{\substack{x\in [G\backslash X]}} | G | [G_x, s(x)]_G = | G | [X,s]\,.
\end{equation*}
The second equality is a similar adaptation of the proof of \cite[Proposition~2.4]{Boltje1998b}.
\end{proof}

\begin{corollary}\label{cor m iso}
Let $(\calG,\calS)$, $\calD$, $F$, and $G$ be as in Definition~\ref{def sigma}.
If $|G|$ is invertible in $R$ then $m_{F,G}$ and $|G|^{-1}n_{F,G}$ are mutually inverse isomorphisms.
\end{corollary}

\begin{corollary}
Let $(\calG,\calS)$, $\calD$, $F$, and $G$ be as in Definition~\ref{def sigma}.
If $F_+(G)$ has trivial $|G|$-torsion then $m_{F,G}$ is injective. 
\end{corollary}


\section{Multiplicative Structures}

\begin{notation}
For any group homomorphism $\varphi\colon H\to K$ we denote by $\lindex{H^\varphi}K_K$ the $(H,K)$-biset $K$ on which $K$ acts by multiplication and $H$ acts via $\varphi$ and multiplication. Similarly, we define the $(K,H)$-biset $\lindex{K}K_{\lexp{\varphi}{H}}$. Note that $\lindex{H^\varphi}K_K\cong (H\times K)/\Dphi$ as $(H,K)$-bisets, where $\Dphi=\{(h,\varphi(h))\mid h\in H\}$, and that $\lindex{K}K_{\lexp{\varphi}{H}}\cong (K\times H)/\phiD$, where $\phiD:=\{(\varphi(h),h)\mid h\in H\}$.
\end{notation}


Green biset functors were defined by Bouc in \cite[\S8.5]{Bouc2010a}. An alternative definition and proof for their equivalence was given by Romero in \cite[Definici\'on~3.2.7, Lema~4.2.3]{RomeroThesis}. In our context, Romero's definition amounts to the following definition.

\begin{definition}\label{def Green biset functor}
Let $(\calG, \calS)$ satisfy axioms (i)--(iii) in \ref{noth GS} and let $\calD = \calC(\calG,\calS)$. 

\smallskip
(a) A {\em Green biset functor} on $\calD$ over $R$ is an object $F\in \calF_{\calD,R}$ together with the datum of an $R$-algebra structure on each $F(G)$, $G\in\calG$, such that the following axioms are satisfied for all $H,K\in\calG$ and all group homomorphisms $\varphi\colon H\to K$:

\smallskip
(i) If $\Dphi \in S(H,K)$ then the map $F([\lindex{K}K_{\lexp{\varphi}{H}}])\colon F(K)\to F(H)$ is a morphism of $R$-algebras.

\smallskip
(ii) If $\Dphi\in S(H,K)$ and
$\phiD\in S(K,H)$ then, for all $a\in F(H)$ and $b\in F(K)$, one has
$$F([\lindex{K}K_{\lexp{\varphi}{H}}])(a)\cdot b = F([\lindex{K}K_{\lexp{\varphi}{H}}])(a\cdot F([\lindex{H^\varphi}K_K])(b))$$
and 
$$b\cdot F([\lindex{K}K_{\lexp{\varphi}{H}}])(a) = F([\lindex{K}K_{\lexp{\varphi}{H}}])
(F([\lindex{H^\varphi}K_K])(b)\cdot a)\,.$$

\smallskip
(b) Let $F_1,F_2\in\calF_{\calD,R}$ be Green biset functors on $\calD$ over $R$. A morphism of Green biset functors between $F_1$ and $F_2$ is a morphism $\eta\colon F_1\to F_2$ in the category $\calF_{\calD,R}$ with the additional property that $\eta_G\colon F_1(G)\to F_2(G)$ is an $R$-algebra homomorphism, for all $G\in\calG$. We denote the category of Green biset functors on $\calD$ over $R$ by $\calF^\mu_{\calD,R}$.
\end{definition}

The following theorem shows that the constructions $F_+$ and $F^+$ also work for Green biset functors. Interestingly, to define the multiplication on $F_+$ one needs restriction maps. The mark morphism is then multiplicative.

\begin{theorem}\label{thm Green +}
Let $(\calG,\calS)$ satisfy the axioms (i)--(iii) in \ref{noth GS} and set $\calD:=\calC(\calG,\calS)$. Further, let $F\in\calF^\mu_{\calD,R}$.

\smallskip
{\rm (a)} Assume that $(\calG,\calS)$  also satisfies axioms (iv) and (vi) in \ref{noth GS} and set $\calD_+:=\calC(\calG,\calS_+)$. Then $F_+\in\calF^\mu_{\calD_+,R}$ and one obtains a functor $-_+\colon \calF^\mu_{\calD,R}\to\calF^\mu_{\calD_+,R}$.

\smallskip
{\rm (b)}  Assume that $(\calG,\calS)$ also satisfies axiom (vii) in \ref{noth GS} and condition $k_2$ in Example~\ref{ex D_+ and D^+}, and set $\calD^+:=\calC(\calG,\calS^+)$. Then $F^+\in\calF^\mu_{\calD^+,R}$ and one obtains a functor $-^+\colon \calF^\mu_{\calD,R}\to\calF^\mu_{\calD^+,R}$.

\smallskip
{\rm (c)} If $(\calG,\calS)$ also satisfies axioms (iv) and (vi) in \ref{noth GS} then the mark morphism $m_{F,G}\colon F_+(G)\to F^+(G)$ is  an $R$-algebra homomorphism for all $G\in\calG$. 
In particular, if $(\calG,\calS)$ additionally satisfies axiom~(vii) in \ref{noth GS} and condition $k_2$ in Example~\ref{ex D_+ and D^+},  we obtain a morphism $m_F\colon F_+\to F^+$ in the category of Green biset functors, $\calF^\mu_{\calE,R}$, where $\calE=\calD^+=\calD_+$.
\end{theorem}

The proof of the previous theorem is straightforward but involves a large number of verifications.
We only point out the multiplicative structures on $F_+(G)$ and $F^+(G)$ and leave the verifications to the reader.
In $F_+(G)$, the product is given by $[X,s]\cdot [Y,t]:=[X\times Y, s\times t]$, where $(X,s),(Y,t)\in\Gamma_F(G)$ and $s\times t\colon  X\times Y\to \coprod_{(x,y)\in X\times Y} F(G_{(x,y)})$ maps $(x,y)$ to
$F(\res_{G_{(x,y)}}^{G_x})(s(x)) \cdot F(\res_{G_{(x,y)}}^{G_y})(t(y))$. This translates to the product 
$$[H,a]_G\cdot[K,b]_G:= \sum_{g\in H\backslash G/K} [H\cap\lexp{x}{K}, F(\res^H_{H\cap\lexp{x}{K}})(a)\cdot
F(\res^{\lexp{x}{K}}_{H\cap\lexp{x}{K}})(\lexp{x}{b})]_G\,,$$
for $H,K\in\Sigma_\calG(G)$ and $a\in F(H)$, $b\in F(K)$. In $F^+(G)$ we define the product coordinate-wise.
In Part~(c), note that the maps $F_+(\res^G_H)$ and $\pi_H\colon F_+(H)\to F(H)$ are $R$-algebra homomorphisms, for all $H\in\Sigma_\calG(G)$.

\begin{example}\label{ex multiplicative local}
(a) Let $G$ be a finite group. Note that the categories $\ConfR(G)$, $\ResfR(G)$, $\MackfR(G)$, considered in Example~\ref{ex D_+ and D^+ local}, have Green biset functor versions $\ConfmuR(G)$, $\ResfmuR(G)$, $\MackfmuR(G)$ which are {\em fused} versions of the $R$-algebra conjugation, restriction and Mackey functors considered in \cite[Section~1]{Boltje1998b}. With this notation, we have again functors $-_+\colon \ResfmuR(G)\to\MackfmuR(G)$ and $-^+\colon \ConfmuR(G)\to\MackfmuR(G)$ that are restrictions of the functors $-_+$ and $-^+$ given in \cite[Section~2]{Boltje1998b}.

\smallskip
(b) The multiplicative structures on the Burnside ring, $A$-fibered Burnside ring, and the global representation ring (see Example~\ref{ex Burnside examples}) coincide with the multiplication on $F_+(G)$ defined above.
\end{example}

Next we turn our attention to {\em species} of $F_+(G)$, with the goal to describe them in terms of species of $F$.
For our purposes, a {\em species} of a ring $\Lambda$ is a ring homomorphism $\sigma\colon \Lambda\to\CC$. We denote the set of species of $\Lambda$ by $\Sp(\Lambda)$.
For the remainder of this section we assume that $R=\ZZ$, that $(\calG,\calS)$ satisfies axioms (i)--(iv) and (vi) in \ref{noth GS}, and that $F\in\calF_{\calD}^\mu:=\calF_{\calD,\ZZ}^\mu$ is a Green biset functor, where $\calD:=\calC(\calG,\calS)$. For $G\in\calG$ we set 
\begin{equation*}
  \Sp_+(F,G):= \{(H,\tau)\mid H\in\Sigma_\calG(G), \tau\in\Sp(F(H))\}
\end{equation*}
and 
\begin{equation*}
  \Sptilde_+(F,G):=\{(H,\tautilde)\mid H\in\Sigma_{\calG}(G), \tautilde\in\Sp(F(H)^{N_G(H)})\}\,,
\end{equation*}
where $F(H)^{N_G(H)}$ denotes the $N_G(H)$-fixed points of $F(H)$, a subring of $F(H)$. 
Note that $G$ acts on $\Sp_+(F,G)$ from the right: For $x\in G$ we define $(H,\tau)^x:=(H^x,\tau^x)$, where $\tau^x:=\tau\circ F(c_x)\colon F(H^x)\to\CC$. Similarly, $G$ acts on $\Sptilde_+(F,G)$ and the map
\begin{equation}\label{eqn rho def}
  \Sp_+(F,G)\to\Sptilde_+(F,G)\,,\quad (H,\tau)\mapsto (H,\tau|_{F(H)^{N_G(H)}})\,,
\end{equation}
is $G$-equivariant. For every $(H,\tautilde)\in\Sptilde_+(F,G)$ we defined the map
\begin{equation*}
  \sigma_{(H,\tautilde)}\colon F_+(G)\to\CC\,,\quad \omega\mapsto\tautilde\bigl( m_{F,G}(\omega)_H\bigr)\,,
\end{equation*}
where $m_{F,G}(\omega)_H$ denotes the $H$-component of $m_{F,G}(\omega)\in(\prod_{H\in\Sigma_\calG(G)}F(H))^G$. Note that $m_{F,G}(\omega)_H\in F(H)^{N_G(H)}$. Since $m_{F,G}$ is a ring homomorphism (see Theorem~\ref{thm Green +}(c)), we have $\sigma_{(H,\tautilde)}\in\Sp(F_+(G))$. It is easy to verify that $\sigma_{(H,\tautilde)}=\sigma_{(H,\tautilde)^x}$ for all $x\in G$. Thus, we obtain a map 
\begin{equation*}
  \sigmatilde_{F,G}\colon \Sptilde_+(F,G)/G\to \Sp(F_+(G))\,, \quad [H,\tautilde]\mapsto \sigmatilde_{(H,\tautilde)}\,,
\end{equation*}
and a commutative diagram
\begin{equation}\label{diag rho sigma}
\parbox{7cm}{
\begin{diagram}
\Sp_+(F,G)/G & \movearrow(10,0){\Ear[30]{\sigma_{F,G}}} & \movevertex(20,0){\Sp(F_+(G))} &&
\Sar{\rho_{F,G}} & &  \movearrow(-10,0){\neaR{\sigmatilde_{F,G}}}  &&
\Sptilde_+(F,G)/G & & &&
\end{diagram}}
\end{equation}
where $\sigma_{F,G}:=\sigmatilde_{F,G}\circ\rho_{F,G}$ and $\rho_{F,G}$ denotes the map induced by the map in (\ref{eqn rho def}). For $(H,\tau)\in\Sp_+(F,G)$ we set $\sigma_{(H,\tau)}:=\sigma_{(H,\tautilde)}$, where $\tautilde:=\tau|_{F(H)^{N_G(H)}}$.

\begin{theorem}
Let $(\calG,\calS)$ satisfy (i)--(iv) and (vi) in \ref{noth GS}, set $\calD:=\calC(\calG,\calS)$, let $F\in\calF_{\calD}^\mu$, and let $G\in\calG$. Then the map $\sigma_{F,G}$ in the commutative diagram (\ref{diag rho sigma}) is injective and the map $\sigmatilde_{F,G}$ is surjective. If moreover $F(H)$ is commutative for all $H\in\Sigma_{\calG}(G)$ then all three maps $\sigma_{F,G}$, $\sigmatilde_{F,G}$ and $\rho_{F,G}$ are bijective.
\end{theorem}

\begin{proof}
We first prove that $\sigma_{F,G}$ is injective.
Let $(H,\tau),(K,\rho)\in \Sp_+(F,G)$ with $\sigma_{(H,\tau)}=\sigma_{(K,\rho)}$. We will show that $(H,\tau)$ and $(K,\rho)$ are $G$-conjugate. Since $\sigma_{(K,\rho)}([H,1]_G) = \sigma_{(H,\tau)}([H,1]_G)  = [N_G(H):H]\neq 0\in \CC$, we have $K\le \lexp{x}{H}$ for some $x\in G$. Similarly, $H\le\lexp{y}{K}$ for some $y\in G$. Thus, $K$ and $H$ are conjugate, and we may assume that they are equal. For every $a\in F(H)^{N_G(H)}$, we have $[N_G(H):H] \tau(a) = s_{(H,\tau)}([H,a]_G) = s_{(K,\rho)}([H,a]_G) = [N_G(H):H] \rho(a)$, which shows that $\tau$ and $\rho$ coincide on $F(H)^{N_G(H)}$. Furthermore, for every $b\in F(H)$, the element $a:=\sum_{x\in [N_G(H)/H]} \lexp{x}{b}$ lies in $F(H)^{N_G(H)}$, so that $\sum_{x\in[N_G(H)/H]}\tau^x(b) = \tau(a) = \rho(a) = \sum_{x\in [N_G(H)/H]} \rho^x(b)$. This implies $\sum_{x\in [N_G(H)/H]}\tau^x=\sum_{x\in [N_G(H)/H]}\rho^x$. Since ring homomorphisms from $F(H)$ to $\CC$ are $\CC$-linearly independent in the $\CC$-vector space of all functions from $F(H)$ to $\CC$, there exists $x\in N_G(H)$ with $\tau=\rho^x$, and the proof of the injectivity of $\sigma_{F,G}$ is complete.

\smallskip
Next we prove that $\sigmatilde_{F,G}$ is surjective. Let $\sigma\in\Sp(F_+(G))$ and let $\sigma_\CC\colon \CC\otimes F_+(G)\to\CC$ be its extension to a $\CC$-algebra homomorphism. By the arguments in the proof of Lemma~7.5(ii) in \cite{Boltje1998b}, one has a canonical $\CC$-algebra isomorphism $\CC\otimes F_+(G)\myiso (\CCF)_+(G)$, where $\CCF\in\calF_{\calD,\CC}^\mu$ denotes the obvious scalar extension of $F$ from $\ZZ$ to $\CC$. We write again $\sigma_\CC\colon (\CCF)_+(G)\to\CC$ for the corresponding $\CC$-algebra homomorphism. Next we consider the commutative diagram
\begin{diagram}
\movevertex(-70,0){F_+(G)} & \movearrow(-65,0){\Ear{m_{F,G}}} & \movevertex(-60,0){F^+(G)} & & 
\movearrow(-30,0){\Ear[140]{\sim}} & & 
\movevertex(30,0){\mathop{\prod}\limits_{H\in\Sigmatilde_\calG(G)} F(H)^{N_G(H)}} &&
\movearrow(-70,0){\sar} & & & & & & \movearrow(40,0){\sar} &&
\movevertex(-70,0){(\CCF)_+(G)} & \movearrow(-65,0){\Ear{m_{\CCF,G}}} & \movevertex(-60,0){(\CCF)^+(G)} & \movearrow(-55,0){\Ear[20]{\sim}} & 
\movevertex(-20,0){ \mathop{\prod}\limits_{H\in\Sigmatilde_{\calG}(G)}(\CC\otimes F(H))^{N_G(H)}} & 
\movearrow(10,0){\War[20]{\sim}} &
\movevertex(50,0){\mathop{\prod}\limits_{H\in\Sigmatilde_\calG(G)}\CC\otimes F(H)^{N_G(H)}} &&
\end{diagram}
where $\Sigmatilde_\calG(G)\subseteq\Sigma_\calG(G)$ denotes a set of representatives of the conjugacy classes of $\Sigma_\calG(G)$ and all morphisms other than $m_{F,G}$ and $m_{\CCF,G}$ are the obvious natural maps, cf.~(\ref{eqn F_+ projection iso}). Note that the right hand map of the bottom row is a $\CC$-algebra isomorphism, since $\CC$ is flat as $\ZZ$-module (see the proof of Lemma~7.5(i) in \cite{Boltje1998b}). Note also that $m_{\CCF,G}$ is a $\CC$-algebra isomorphism by Theorem~\ref{thm Green +}(c) and Corollary~\ref{cor m iso}. Thus, all maps in the bottom row are $\CC$-algebra isomorphisms and it follows that $\sigma_\CC\colon (\CCF)_+(G)\to \CC$ must come via composition with these isomorphisms from a $\CC$-algebra homomorphism $\CC\otimes F(H)^{N_G(H)}\to \CC$, for some $H\in\Sigmatilde_{\calG}(G)$. The restriction of this homomorphism to $F(H)^{N_G(H)}$ yields a species $\tautilde\in\Sp(F(H)^{N_G(H)})$ with $\sigma=\sigma_{(H,\tautilde)}$. Thus, $\sigmatilde_{F,G}$ is surjective.

\smallskip
Assume for the rest of the proof that $F(H)$ is commutative for all $H\in\Sigma_\calG(G)$. By Theorem~1.8.1 and the first part of the proof of Theorem 2.9.1 in \cite{Benson1984}, every $\CC$-algebra homomorphism $\tautilde\colon (\CC\otimes F(H))^{N_G(H)}\to\CC$ extends to a $\CC$-algebra homomorphism $\tau\colon \CC\otimes F(H)\to\CC$. Since $(\CC\otimes F(H))^{N_G(H)}\cong \CC\otimes F(H)^{N_G(H)}$ as $\CC$-algebras (with the same argument as above), this implies that the map in (\ref{eqn rho def}) and therefore also the map $\rho$ in Diagram~(\ref{diag rho sigma}) is surjective. Since $\sigma_{F,G}$ is injective, also $\rho_{F,G}$ is injective. Thus, $\rho$ is bijective. Since $\sigma_{F,G}$ is injective, $\sigmatilde_{F,G}$ is surjective and $\rho_{F,G}$ is bijective, all three maps must be bijective. This completes the proof of the Theorem.
\end{proof}

\section{Adjointness}

\begin{notation}\label{not adj D}
Unless otherwise stated,
throughout this section let $R$ denote a commutative ring.
Let $(\calG,\calS)$ satisfy axioms (i) -- (iv), and set $\calD:=\calC(\calG,\calS)$.
For $H,K\in\calG$ set $\calS_-(H,K) = \{ D \in \calS(H,K) \mid p_1(D) = H \}$.
\end{notation}

The proof of the following lemma is straightforward and is left to the reader.

\begin{lemma}
Let $(\calG,\calS)$ be as in \ref{not adj D}. Then $(\calG,\calS_-)$ also satisfy axioms (i)--(iv) and
$(\calS_-)_+ = \calS_+$. We can thus define $\calD_- = \calC(\calG,\calS_-)$ and obtain 
$\calD_-\subseteq\calD\subseteq\calD_+$. For any $F\in\calF_{\calD_-,R}$ one has a morphism
$\eta_F\colon F\to \Res^{\calD_+}_{\calD_-}(F_+)$ in $\calF_{\calD_-,R}$ given by 
$\eta_{F,G}(a):=[G,a]_G$ for any $G\in\calG$ and $a\in F(G)$.
\end{lemma}

\begin{theorem}\label{thm adjunction}
Let $(\calG,\calS)$ be as in \ref{not adj D}. Then the functor $-_+\colon \calF_{\calD_-,R}\to \calF_{\calD_+,R}$ is left adjoint to the restriction functor $\Res^{\calD_+}_{\calD_-}\colon \calF_{\calD_+,R}\to\calF_{\calD_-,R}$. More precisely, for any $F\in\calF_{\calD_-,R}$ and $M\in\calF_{\calD_+,R}$, the map $\varphi\mapsto \varphi\circ\eta_F$ defines an $R$-linear isomorphism
\begin{equation*}
  \Hom_{\calF_{\calD_+,R}}(F_+,M)\myiso \Hom_{\calF_{\calD_-,R}}(F,\Res^{\calD_+}_{\calD_-}(M))\,.
\end{equation*}
%
%
\end{theorem}

\begin{proof}
We will show that for any $\psi\in\Hom_{\calF_{\calD_-,R}}(F,\Res^{\calD_+}_{\calD_-}(M))$ there exists a unique 
$\varphi\in  \Hom_{\calF_{\calD_+,R}}(F_+,M)$ with $\psi=\varphi\circ\eta_F$. We show the uniqueness first. Let $G\in\calG$. Then every element in $F_+(G)$ can be written as an $R$-linear combination of elements of the form $[H,a]_G=F_+(\ind_H^G)([H,a]_H)$ with $a\in F(H)$ and $H\in\Sigma_{\calG}(G)$. Since $[H,a]_H$ is in the image of $\eta_{F,H}$ and $\varphi$ commutes with inductions, $\varphi$ is uniquely determined by the condition $\psi=\varphi\circ\eta_F$. Next we show the existence of $\varphi$. For $G\in\calG$ and an object $(X,s)$ in $\Gamma_F(G)$, we define $\varphi_G([X,s]):=\sum_{x\in [G\backslash X]} M(\ind^G_{G_x}) (\psi_{G_x}(s(x)))$. Note that this yields a well-defined map $\varphi_G\colon F_+(G)\to M(G)$. In fact, the above sum does not depend on the choice of $[G\backslash X]$, and by the definition of $F_+(G)$ we only have to check that the relations in Definition~\ref{def 2 operations} are respected, which is an easy verification. Note also that, for $H\in\Sigma_\calG(G)$ and $a\in F(H)$, we have $\varphi_G([H,a]_G)=M(\ind_H^G)(\psi_H(a))$. Choosing $H=G$, this shows that $\varphi\circ\eta_F=\psi$. Next we show that $\varphi\in\calF_{\calD_+,R}$. So let $G,H\in\calG$, $D\in\calS_+(G,H)$ and set $U:=(G\times H)/D$. Let $L\in\Sigma_\calG(H)$ and $a\in F(L)$. On the one hand we have
\begin{eqnarray*}
M([U])(\varphi_H([L,a]_H)) & = & M([U])\left(\ind^H_{L} (\psi_{L}(a)) \right)\\
  & = &   M\left(\left[ \frac{G\times H}{D} \times_H \frac{H\times L}{\Delta(L)}  \right] \right) (\psi_L(a))\\
& = &  \sum_{x\in [p_2(D)\backslash H / L]} M\left( \left[\frac{G\times L}{D*\lexp{(x,1)}{\Delta(L)}} \right]\right)(\psi_L(a))
\end{eqnarray*}
On the other hand, using (\ref{eqn explicit F_+}), we have
\begin{equation*}
  F_+([U])([L,a]_H) = \sum_{x\in[p_2(D)\backslash H/L]} 
  \left[ D*\xL, F\left(\left[\frac{D*\xL \times \xL}{D*\Delta(\xL)}\right]\right)(\lexp{h}{a})\right]_G
\end{equation*}
and therefore
\begin{align*}
 \varphi_G\bigl(F_+([U])([L,a]_H)\bigr) & = \sum_{x\in[p_2(D)\backslash H/L]}
     M(\ind_{D*\xL}^G)\bigl(\psi_{D*\xL}\bigl(F\left(\left[\frac{D*\xL\times\xL}{D*\Delta(\xL)}\right]\right)(\lexp{x}{a})\bigr)\bigr)\\
  & = \sum_{x\in[p_2(D)\backslash H/L]}
     M(\ind_{D*\xL}^G)\bigl( M\left(\left[\frac{D*\xL\times\xL}{D*\Delta(\xL)}\right]\right) (\psi_{\xL}(\lexp{x}{a}))\bigr)\\
  & = \sum_{x\in[p_2(D)\backslash H/L]}
        M\left(\left[\frac{G\times D*\xL}{\Delta(D*\xL)}\mathop{\cdot}\limits_{D*\xL}
           \frac{D*\xL\times\xL}{D*\Delta(\xL)}\right]\right)(\psi_{\xL}(\lexp{x}{a}))\\
  & = \sum_{x\in[p_2(D)\backslash H/L]}
        M\left(\left[\frac{G\times \xL}{\Delta(D*\xL)*D*\Delta(\xL)}\right]\right)(\psi_{\xL}(\lexp{x}{a}))\,.
\end{align*}
Since $\psi_{\xL}(\lexp{x}{a}) = M(c_x)(\psi_L(a))$ and $\Delta(D*\xL)*D*\Delta(\xL)*\lexp{(x,1)}{\Delta(L)} = D*\lexp{(x,1)}{\Delta(L)}$, we obtain $M([U])(\varphi_H([L,a]_H)) = \varphi_G(F_+([U])([L,a]_H))$. Thus, $\varphi\in
\Hom_{\calF_{\calD_+,R}}(F_+,M)$, and the proof is complete.
\end{proof}

If $F\in\calF_{\calD_-,R}^\mu$ is a Green biset functor then the natural transformation $\eta_F$ is multiplicative, 
i.e., a morphism in $\Hom_{\calF_{\calD_-,R}^\mu}(F,\Res^{\calD^+}_{\calD_-}(F^+))$. 
Theorem~\ref{thm adjunction} has the following multiplicative version.

\begin{theorem}\label{thm adjunction mu}
Assume that $(\calG,\calS)$ satisfies Axioms (i)--(iv) and Axiom (vi) in \ref{noth GS}. Then the functor $-_+\colon \calF^\mu_{\calD_-,R}\to \calF^\mu_{\calD_+,R}$ is left adjoint to the restriction functor $\Res^{\calD_+}_{\calD_-}\colon \calF^\mu_{\calD_+,R}\to\calF^\mu_{\calD_-,R}$. As in Theorem~\ref{thm adjunction}, the adjunction bijection is given by composition with $\eta_F$.
\end{theorem}

\begin{proof}
Let $F\in\calF^\mu_{\calD_-,R}$ and $M\in\calF^\mu_{\calD_+,R}$ and let $\psi\in\Hom_{\calF^\mu_{\calD_-,R}}(F,M)$. We define $\varphi\in\Hom_{\calF_{\calD_+,R}}(F_+,M)$ as in the proof of Theorem~\ref{thm adjunction}. It suffices to show that, for any $G\in\calG$, the map $\varphi_G\colon F_+(G)\to M(G)$ is multiplicative. Let $K,L\in\Sigma_{\calG}(G)$, $a\in F(K)$, and $b\in F(L)$. On the one hand we have
\begin{align*}
  & \ \ \varphi_G([K,a]_G\cdot[L,b]_G) \\
   = &  \sum_{\substack{x\in [ K\backslash G/L ]}} \varphi_G([K\cap\xL, F(\res^{K}_{K\cap \xL})(a) \cdot
  \res^{\xL}_{K\cap\xL}(\lexp{x}{b})]_G) \\
   = &   \sum_{\substack{x\in  [ K\backslash G/L ]}} M(\ind^G_{K\cap\xL})
  \left(\psi_{K\cap\xL}( F(\res^{K}_{K\cap\xL})(a)) \cdot
 \psi_{K\cap\xL}( F(\res^{\xL}_{K\cap \xL})(\lexp{x}{b}))\right)\,.
 \end{align*}
On the other hand
\begin{align*}
  & \ \ \varphi_G([K,a]_G)\cdot \varphi_G([L,b]_G)  =  M(\ind^G_K) (\psi_K(a)) \cdot M(\ind^G_L) (\psi_L(b))\\
  = & \ \ M(\ind^G_K)\left(\psi_K(a)\cdot M(\res^G_K  \mathop{\cdot}\limits_G \ind^G_L)(\psi_L(b)) \right)\\
  =  &   \sum_{\substack{x\in [ K\backslash G/L ]}} M(\ind^G_{K})
  \left(\psi_{K}(a)\cdot M(\ind^K_{K\cap\xL}) \bigl( M(\res^{\xL}_{K\cap\xL})
  (\psi_{\xL}(\lexp{x}{b})) \bigr) \right) \\
  = & \sum_{\substack{x\in [ K\backslash G/L ]}} M(\ind^G_{K})
   \bigl(M(\ind^K_{K\cap\xL})\bigl( M(\res^K_{K\cap\xL}) ( \psi_K(a))\cdot 
   M(\res^{\xL}_{K\cap\xL})(\psi_{\xL}(\lexp{x}{b}))\bigr)\bigr)\,.
 \end{align*}
Since induction is transitive and $\psi$ commutes with restrictions (note that $\calD_-$ contains all possible restrictions), the two expressions coincide.
\end{proof}



\end{document}